\documentclass[10pt,reqno]{amsart}
\usepackage{amsaddr}
\usepackage{graphicx}
\usepackage{amssymb}
\usepackage{amsfonts}
\usepackage[]{amsmath}
\usepackage[]{epsfig}
\usepackage[]{pstricks}
\newpsobject{malla}{psgrid}{subgriddiv=1,griddots=10,gridlabels=6pt}
\usepackage[]{float}
\usepackage{setspace}
\usepackage{tikz}
\usepackage{listings}
\usepackage{comment}
\usepackage{subfigure}
\usepackage{lipsum}
\usepackage[backref=page]{hyperref}
\usepackage{footnote}
\usepackage{orcidlink}
\usetikzlibrary{decorations.shapes}
\newpsobject{malla}{psgrid}{subgriddiv=1,griddots=10,gridlabels=6pt}
\newtheorem{theorem}{Theorem}[section]
\newtheorem{lemma}{Lemma}[section]

\newtheorem{example}{Example}[section]
\numberwithin{equation}{section}
\theoremstyle{remark}
\newtheorem{remark}{Remark}[section]
\newcommand{\px}{\partial_x}
\newcommand{\pt}{\partial_t}
\newcommand{\R}{\mathbb R}

\newcommand\comentario[1]\null
\newlength\replength
\newcommand\repfrac{.33}

\newcommand\rulewidth{.6pt}
\newcommand\tdashfill[1][\repfrac]{\cleaders\hbox to \replength{
  \smash{\rule[\arraystretch\ht\strutbox]{\repfrac\replength}{\rulewidth}}}\hfill}

\newcommand\tdotfill[1][\repfrac]{\cleaders\hbox to \replength{
\smash{\raisebox{\arraystretch\dimexpr\ht\strutbox-.1ex\relax}{.}}}\hfill}

\usepackage{booktabs}
\usepackage{arydshln}

\def\norm#1{\|#1\|}
\def\bra#1{\langle#1\rangle}
\def\wt#1{\widetilde{#1}}

\def\set#1{\{#1\}}
\AtEndDocument{\noindent Márcio Cavalcante and José Marques\\
Instituto de Matemática\\
Universidade Federal de Alagoas\\
Maceió\\
Brazil\\
e-mail: marcio.melo@im.ufal.br\\

\noindent Chulkwang Kwak\\
Department of Mathematics\\
Ewha Womans University\\
Seoul\\
Korea\\
e-mail: ckkwak@ewha.ac.kr\\

\noindent José Marques\\
e-mail: jose.neto@im.ufal.br}
\begin{document}
\pagenumbering{arabic}	
\title[Kawahara on star graphs]{The Kawahara equation  on star graphs}
\author[Cavalcante]{M\'arcio Cavalcante \orcidlink{0000-0001-9873-1765}}
\address{\emph{Instituto de Matem\'atica, Universidade Federal de Alagoas,\\ Macei\'o-Brazil}}
\author[Kwak]{Chulkwang Kwak \orcidlink{0000-0001-9634-6548}}
\address{\emph{Department of Mathematics, Ewha Womans University,\\ Seoul-Korea}}
\let\thefootnote\relax
\author[Marques]{José Marques \orcidlink{0009-0001-9009-2628}}
\address{\emph{Instituto de Matem\'atica, Universidade Federal de Alagoas,\\ Macei\'o-Brazil}}

\begin{abstract}
In this paper, we establish local well-posedness for the Cauchy problem associated with the Kawahara equation on a general metric star graph. Initially, we identify suitable boundary conditions that produce a well-behaved dynamics for the linear equation. Subsequently, we derive the integral formula using the forcing operator method, previously applied to the Kawahara equation on the half-line by Cavalcante and Kwak \cite{Cavalcantekwak2020}, and the Fourier restriction method of Bourgain \cite{Bourgain1}. This work has the potential to be extended to other fifth-order nonlinear dispersive equations on star graphs.
\medskip

\end{abstract}
\maketitle
\section{Introduction}

Over the past few decades, models based on partial differential equations have been very effective in tackling many problems dealing with flows on networks (e.g., irrigation channels, gas pipelines, blood circulation, vehicular traffic, supply chains, air traffic management – see \cite{berkolaiko,Bressan} for  surveys of the topic).

In the context of nonlinear dispersive equation on star graphs,  some results have appeared in the literature in the last years. For the Korteweg-de Vries (KdV) equation on star graphs, see the work \cite{CM} and the references therein. For the Biharmonic equation, we refer the reader to \cite{fourth}. For the Nonlinear Schrödinger (NLS) equation, see the works \cite{Adami,Angulo,Ardila,Nojanls} and references therein. For the Benjamin–Bona–Mahony (BBM) equation, see Bona and Cascavel \cite{bbm}. Also, in \cite{Caudrelier}, a framework was presented to solve the open problem of formulating the inverse scattering method (ISM) for an integrable PDE on a star-graph, and the nonlinear Schrödinger equation was used to illustrate the method. More recently, the small amplitude linearization of the Korteweg-de Vries equation on the line with a local defect scattering waves, represented by a metric graph domain adjoined at one point, was studied in \cite{smith}, using the unified transform method on a metric graph. 

In this work, we consider the following Kawahara equation:
\begin{equation}\label{eq:5kdv}
\pt u - \px^5 u + \px(u^2)=0.
\end{equation}
The Kawahara equation was first proposed by Kawahara \cite{Kawahara1972} describing solitary-wave propagation in media. Also, the Kawahara equation can be described in the theory of magneto-acoustic sound waves in plasma and in theory of shallow water waves with surface tension. For another physical background of Kawahara equation or in view of perturbed equation of KdV equation, see \cite{HS1988, PRG1988, Boyd1991} and references therein.

The Cauchy problem for the Kawahara equation on $\R$ has been extensively studied on the works \cite{CDT2006,CT2005,WCD2007,CLMW2009,JH2009,CG2011,Kato2011,Himonas} and references therein. 

The literature also addresses problems defined on positive and negative half-lines. Specifically, initial boundary value problems (IBVPs) for the Kawahara equation on the right, given by:
\begin{equation}\label{kawahararight}
\begin{cases}
\pt u - \px^5 u + \px(u^2)=0, & (t,x)\in (0,T)  \times(0,\infty),\\
u(0,x)=u_0(x),                                   & x\in(0,\infty),\\
u(t,0)=f(t),\ u_x(t,0)=g(t)& t\in(0,T)
\end{cases}
\end{equation}
 can be seen in the works \cite{San2003, Lar1,Lar2,Cavalcantekwak2020}.

The IBVP on the left half-line $(-\infty, 0)$
\begin{equation}\label{kawaharaleft}
\begin{cases}
	\pt u - \px^5 u + \px(u^2)=0, & (t,x)\in  (0,T) \times (-\infty,0),\\
	u(0,x)=u_0(x),                                   & x\in(-\infty,0),\\
	u(t,0)=f(t),\ 	u_x(t,0)=g(t),\ u_{xx}(t,0)=h(t)& t\in(0,T).
\end{cases}
\end{equation}
was considered in \cite{Cavalcantekwak2020}.

In particular, in \cite{Cavalcantekwak2020} local well-posedness was obtained under low regularity assumptions, by using the approach introduced in \cite{CK}, which is based on an extension problem using Riemann-Liouville fractional integration. As far as we know, the Kawahara  equation on star graphs has not been previously explored.

 \subsection{Formulation of the problem}

The present paper focuses on the study of the Kawahara equation on general star graphs, specifically addressing the scenario where all edges have unbounded lengths. This research investigates a problem that, to the best of our knowledge, has not been previously addressed in the literature.

In this spirit, we consider the Kawahara equation on a general star graph given by $k$ negative half-lines and $m$ positive half-lines attached to a common point, called the vertex. We denote this graph by $\mathcal G$, composed of $k$ negative half-lines and $m$ positive half-lines, attached by a unique vertex. The vertex corresponds to $0$. We consider the following system of Kawahara equations
\begin{equation}\label{Kawaharaprevia}
\begin{cases}
\partial_t u_i-\partial^5_{x}u_i+\partial_x(u_i^2)=0,& (x, t)\in L_i\times (0,T),\ i=1,2, ..., k\\
\partial_t v_j-\partial^5_{x}v_j+\partial_x(v_j^2)=0,& (x,t)\in R_j\times (0,T),\ j=1,2, ..., m.\\
\end{cases}
\end{equation}
We will need to impose $3k+2m$ boundary conditions, which should also provide connection between the bonds (see next subsection).
Furthermore, we assume the initial conditions
\begin{align}\label{ICprevia}
u_i(x,0)&=u_{0i},\ x\in L_i,\ i=1,2, ..., k,\nonumber\\
v_j(x,0)&=v_{0j},\ x\in R_j,\ j=1,2, ..., m.
\end{align}

\begin{figure}[htp]\label{stark+m}
	\centering 
	\begin{tikzpicture}[scale=2]
\node at (0,1)[rotate=90]{\Huge$\vdots$};
\node at (0,-1)[rotate=90]{\Huge$\vdots$};
	
\foreach \angle in {30,60,120,150} {
    \draw[dashed, red] (\angle:2cm) -- (0,0);
    \draw[->,thick,red] (\angle:0cm) -- (\angle:1.5cm);
}

\foreach \angle in {210,240,300,330} {
    \draw[dashed, black] (\angle:2cm) -- (0,0);
    \draw[->,thick,black] (\angle:1.5cm) -- (\angle:1.4cm);
    \draw[-,thick, black] (\angle:1.5cm) -- (\angle:0.0cm);
}
	\end{tikzpicture}
	\caption{A star graph with $m+k$ edges}
	\end{figure}
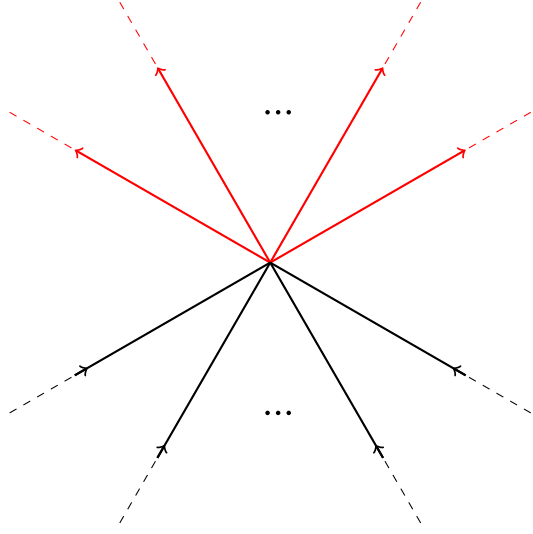

\subsection{Choices of boundary conditions} The selection of appropriate boundary conditions for a well-posed problem in star graphs, in general, is a complex task. To the best of our knowledge, different from the Schrödinger equation, the physically relevant boundary conditions for the Kawahara equation on a star graph remain an open question. Therefore, in this study, we focus on a class of boundary conditions which are consistent with the linearized Kawahara equation. Equivalently, the boundary conditions chosen here ensure uniqueness for smooth, decaying solutions of the linearized Cauchy problem \eqref{Kawaharaprevia}-\eqref{ICprevia}.

Here, we consider star graphs $\mathcal{G},$ composed by $k+m$ edges in a set of edges $E=E_{-}\cup E_{+}$, where $E_{-}=\{L_i\}_{i=1}^k$ and $E_{-}=\{R_i\}_{i=1}^m$ ($L_i=(-\infty,0),\ R_i=(0,\infty)$).

In this spirit, we perform formal computations to obtain appropriate boundary conditions. Initially, let  ${u}(x,t)=(u_1(x,t), ..., u_{k}(x,t))$ and ${v}(x,t)=(v_1(x,t), ..., v_m(x,t)),$ and suppose that $(\overrightarrow{u},\overrightarrow{v})$ is a smooth decaying solution   of the linear version of \eqref{Kawaharaprevia}-\eqref{ICprevia} with homogeneous initial condition, i.e.
\begin{equation}\label{Kawaharalinear}
\begin{cases}
\partial_t u_i-\partial^5_{x}u_i=0,& (x, t)\in L_i\times (0,T),\ i=1,2, ..., k;\\
\partial_t v_j-\partial^5_{x}v_j=0,& (x,t)\in R_j\times (0,T),\ j=1,2, ..., m;\\
u_i(x,0)=u_{0i}(x),&  x \in  L_i, (i=1,2,..., k;\\
v_j(x,0)=v_{0j}(x),&  x \in  R_i, (i=1,2,..., m).
\end{cases}
\end{equation}

We will prove that the corresponding solution to the problem $\eqref{Kawaharalinear}$ is identically zero, which implies the uniqueness argument. To this end, multiplying the first equation in \eqref{Kawaharalinear} by $u,$ the second by $v$, integrating by parts and then using the initial condition, we infer that
\begin{equation}\label{uniqueness}
\begin{split}
&\sum_{i=1}^{k} \int_{-\infty}^{0} (u_i)^2(x,T) dx + \sum_{j=1}^{m} \int_{0}^{\infty} (v_j)^2(x,T) dx=\\
&\sum_{i=1}^{k} \int_{0}^{T} \left( \partial^2_{x} u_i \right)^2 (0,t) dt - \sum_{j=1}^{m} \int_{0}^{T} (\partial^2_{x} v_j)^2 (0,t)) dt\\
&- 2 \sum_{i=1}^{k} \int_{0}^{T} u_i(0,t) \partial^4_{x} u_i(0,t) dt + 2 \sum_{j=1}^{m} \int_{0}^{T} v_j(0,t) \partial^4_{x} v_j(0,t) dt\\
&+ 2 \sum_{i=1}^{k}\int_0^T \partial_{x} u_i (0,t) \partial^3_{x} u_i(0,t) dt - 2 \sum_{j=1}^{m}\int_0^T \partial_{x} v_j (0,t) \partial^3_{x} v_j(0,t) dt.
\end{split}
\end{equation}

By analyzing \eqref{uniqueness}, we are interested in boundary conditions for the Cauchy problem \eqref{Kawaharaprevia}-\eqref{ICprevia} such that the right-hand side of \eqref{uniqueness} would have a nonpositive sign. In this sense, if we consider the following particular boundary conditions
\begin{align}
&u_1(0,t)=a_i u_i(0,t)=a_{k+j} v_j(0,t),\ &i=2, ..., k,\ j=1, ..., m \label{a}\\
&\partial_x u_1(0,t)=b_i \partial_x u_i(0,t)=b_{k+j} \partial_x v_j(0,t),\ &i=2, ..., k,\ j=1, ..., m \label{b}\\
&\partial^2_x U(0,t)=C\partial^2_x V(0,t),& \label{c}\\
&\displaystyle\sum^k_{i=1} b_i^{-1}\partial^3_{x} u_i(0,t)=\displaystyle\sum^m_{j=1} b^{-1}_{k+j} \partial^3_{x} v_j(0,t),& \label{d}\\
&\displaystyle\sum^k_{i=1} a_i^{-1}\partial^4_{x} u_i(0,t)=\displaystyle\sum^m_{j=1} a^{-1}_{k+j} \partial^4_{x} v_j(0,t),& \label{e}
\end{align}
where   $V(x,t)=(v_1(x,t), ..., v_{m}(x,t))^T,\ U(x,t)=(u_1(x,t), ..., u_{k}(x,t))^T$ and $C$ is a $k \times m$ matrix such that {($C^TC-I_m$)} is a non-positively defined matrix (equivalently, $C$ is a contraction, i.e., $\| C(x)\|\leq\|x\|\ \forall x\in \R^m$). $a_i,\ a_{k+j},\ b_i,\ b_{k+j}$ are non-zero constants. Without loss of generality, we will consider $a_1=b_1=1$.

From  conditions \eqref{a}-\eqref{c}, we get that
\begin{equation}\label{uniqueness1}
\begin{split}\\
&\sum_{i=1}^{k} \int_{-\infty}^{0} (u_i)^2(x,T) dx + \sum_{j=1}^{m} \int_{0}^{\infty} (v_j)^2(x,T) dx\\
&= \int_0^T (\partial^2_{x} V(0,t))^T (C^T C - I_m) \partial^2_{x} V(0,t) dt\\
&- 2 \sum_{i=1}^k \int_0^T a_i^{-1} u_1(0,t) \partial^4_{x} u_i(0,t) dt + 2 \sum_{j=1}^m \int_0^T a_{k+j}^{-1} u_1(0,t) \partial^4_{x} v_j(0,t) dt\\
&+ 2 \sum_{i=1}^k \int_0^T b_i^{-1} \partial_{x} u_1(0,t) \partial^3_{x} u_i(0,t) dt - 2 \sum_{j=1}^m \int_0^T b_{k+j}^{-1} \partial_{x} u_1(0,t) \partial^3_{x} v_j(0,t) dt.
\end{split}
\end{equation}

Now, from identities \eqref{d}-\eqref{e} and the hypothesis over $C$ we obtain that
\begin{equation}\label{uniqueness2}
\begin{split}
\sum_{i=1}^k\int_{-\infty}^0(u_i)^2(x,T)dx+\sum_{j=1}^{m}\int_{0}^{+\infty}(v_j)^2(x,T)dx=
\lambda \sum_{j=1}^m\int_0^T(\partial^2_{x}v_j)^2(0,t)dt,
\end{split}
\end{equation}

for some $\lambda\leq 0$. It follows that ${u}(x,T)={0}$ and ${v}(x,T)={0}$.

\subsection{Setting of the problem}
Based in the previous discussion, we consider the following Cauchy problem associated to Kawahara equation on star graphs $\mathcal{G}$ given by
\begin{equation}\label{Kawahara}
\begin{cases}
\partial_t u_i-\partial^5_{x}u_i+\partial_x(u_i^2)=0,& (x, t)\in L_i\times (0,T),\ i=1,2, ..., k\\
\partial_t v_j-\partial^5_{x}v_j+\partial_x(v_j^2)=0,& (x,t)\in R_j\times (0,T),\ j=1,2, ..., m.\\
u_i(x,0)=u_{0i}(x),&  x \in  L_i, (i=1,2,..., k)\\
v_j(x,0)=v_{0j}(x),&  x \in  R_i, (i=1,2,..., m),\\
\end{cases}
\end{equation}
with the following set of boundary conditions.
\begin{equation}\label{condfront}
\begin{cases}
u_1 = a_i u_i = a_{k+j} v_j \quad (i=2,\hdots, k) (j=1,\hdots,m);  & t\in (0,T)\\
    \partial_x u_{1} = b_i \partial_x u_{i} = b_{k+j} \partial_x v_{j} \quad (i=2,\hdots, k) (j=1,\hdots,m);  & t\in (0,T)\\
    \partial^2_{x} U = C \partial^2_{x} V; & t\in (0,T)\\
    \sum_{i=1}^k b_i^{-1} \partial^3_{x} u_i = \sum_{j=1}^m b_{k+j} \partial^3_{x} v_j; & t\in (0,T)\\
    \sum_{i=1}^k a_i^{-1} \partial^4_{x} u_i = \sum_{j=1}^m a_{k+j}^{-1} \partial^4_{x} v_j & t\in (0,T),
\end{cases}
\end{equation}
where   $V(x,t)=(v_1(x,t), ..., v_{m}(x,t))^T,\ U(x,t)=(u_1(x,t), ..., u_{k}(x,t))^T$ and $C=[c_{ij}]_{k\times m}$ is a $k \times m$ matrix such that  $$\| C(x)\|\leq\|x\|\ \forall x\in \R^m,$$  and  $a_i,\ a_{k+j},\ b_i,\ b_{k+j}$ are non-zero constants. Without loss of generality, we will consider $a_1=b_1=1$.

Here, we consider the initial data on the Sobolev spaces, i.e., $$(U_0(x),V_0(x)):=((u_{0,1}, u_{0,2},...,u_{0,k}),(v_{0,1}, v_{0,2},...,v_{0,k}))\in H^{s}(\mathcal{G})=\bigoplus_{e \in E_-} H^s(-\infty, 0) \oplus \bigoplus_{e \in E_+} H^s(0, +\infty).$$

Our goal in studying the Cauchy problem \eqref{Kawahara}-\eqref{condfront} is to obtain results of local well-posedness, in Sobolev spaces with low regularity.

Initially, as mentioned previously, we denote the classical $L^2$ based Sobolev spaces on $\mathcal G$ by
\[ H^s(\mathcal{G}) = \bigoplus_{e \in E_-} H^s(-\infty, 0) \oplus \bigoplus_{e \in E_+} H^s(0, +\infty) \]

Here, we are interested in solving the problem in the set of regularity $s\in (-\frac12,\frac 52)\setminus\{ \frac12,\frac32\}$. Also, for $s>\frac12$, we consider the following compatibility condition
\begin{equation}\label{compatibility1}
    u_{0,1} (0)= a_i u_{0,i} (0)= a_{k+j} v_{0,j} (0);
\end{equation}
and, for $s>\frac32$, this other condition
\begin{equation}\label{compatibility2}
    \partial_x u_{0,1} (0)= b_i \partial_x u_{0,i} (0)= b_{k+j} \partial_x v_{0,j} (0).
\end{equation}

\section{Main result}
To introduce the main theorem, we shall make the notation clear. In this sense, in this section, the indices $i,j,\vartheta,l$ and $p$ are integers. $i$ is between $1$ and $k$; $j$ is between $1$ and $m$; $\vartheta$ is between $0$ and $4$; $l$ is between $1$ and $3$; and $p$ is $1$ or $2$.
\begin{equation}\label{genmatrix}
\mathbf{M}(\boldsymbol\lambda, \boldsymbol\beta)=\left[\begin{array}{cc}
M_{1,1} & 0\\
M_{2,1} & M_{2,2}\\
M_{3,1} & 0\\
M_{41} & M_{42}\\
M_{51} & M_{52}\\
M_{61} & M_{62}\\
M_{71} & M_{72}\\
\end{array}\right]
\end{equation}

\noindent is a $(3k+2m) \times (3k+2m)$ matrix formed by the following blocks, whose orders are shown subscript,\\[15pt]
\begin{equation}\label{m11}
M_{1,1}=\left[\begin{array}{ccccccccccccccccccccccccccc}
\rho^{(0)}_{1,1} & \rho^{(0)}_{1,2} & \rho^{(0)}_{1,3} & -a_{2}\rho^{(0)}_{2,1} & -a_{2}\rho^{(0)}_{2,2} & -a_2\rho^{(0)}_{2,3} & & 0 & 0 & 0\\
    \vdots & \vdots & \vdots & & & &\ddots & & & \\
    \rho^{(0)}_{1,1} & \rho^{(0)}_{1,2} & \rho^{(0)}_{1,3} & 0 & 0 & 0 & & -a_{k}\rho^{(0)}_{k1} & -a_{k}\rho^{(0)}_{k2} & -a_k\rho^{(0)}_{k3}
\end{array}\right]_{(k-1) \times 3k},\\[5pt]
\end{equation}
\begin{equation}\label{m21}
M_{2,1}=\left[\begin{array}{cccccccccccccccccccccccccccccc}
\rho^{(0)}_{1,1} & \rho^{(0)}_{1,2} & \rho^{(0)}_{1,3} & 0 & \hdots & 0\\
    \vdots & \vdots & \vdots & \vdots & \ddots & \vdots\\
    \rho^{(0)}_{1,1} & \rho^{(0)}_{1,2} & \rho^{(0)}_{1,3} & 0 & \hdots & 0\\
\end{array}\right]_{m\times 3k},\\[5pt]
\end{equation}
\begin{equation}\label{m22}
M_{2,2}=\left[\begin{array}{ccccccccccccccccccccccccccc}
-a_{k+1}\iota^{(0)}_{1,1} & -a_{k+1}\iota^{(0)}_{1,2} & & 0 & 0\\
& & \ddots & & \\
0 & 0 &  & -a_{k+m}\iota^{(0)}_{m1} & -a_{k+m}\iota^{(0)}_{m2}
\end{array}\right]_{m \times 2m},
\end{equation}

\begin{equation}\label{m31}
M_{3,1}=\left[\begin{array}{ccccccccccccccccccccccccccc}
\rho^{(1)}_{1,1} & \rho^{(1)}_{1,2} & \rho^{(1)}_{1,3} & -b_{2}\rho^{(1)}_{2,1} & -b_{2}\rho^{(1)}_{2,2} & -b_2\rho^{(1)}_{2,3} & & 0 & 0 & 0\\
    \vdots & \vdots & \vdots & & & &\ddots & & &\\
    \rho^{(1)}_{1,1} & \rho^{(1)}_{1,2} & \rho^{(1)}_{1,3} & 0 & 0 & 0 & & -b_{k}\rho^{(1)}_{k1} & -b_{k}\rho^{(1)}_{k2} & -b_k\rho^{(1)}_{k3}
\end{array}\right]_{(k-1) \times 3k},\\[5pt]
\end{equation}

\begin{equation}\label{m41}
M_{4,1}=\left[\begin{array}{cccccccccccccccccccccccccccccc}
\rho^{(1)}_{1,1} & \rho^{(1)}_{1,2} & \rho^{(1)}_{1,3} & 0 & \hdots & 0\\
    \vdots & \vdots  & \vdots  & \vdots  & \ddots & \vdots \\
    \rho^{(1)}_{1,1} & \rho^{(1)}_{1,2} & \rho^{(1)}_{1,3} & 0 & \hdots & 0\\
\end{array}\right]_{m\times 3k},\\[5pt]
\end{equation}
\begin{equation}\label{m42}
M_{4,2}=\left[\begin{array}{ccccccccccccccccccccccccccc}
-b_{k+1}\iota^{(1)}_{1,1} & -b_{k+1}\iota^{(1)}_{1,2} & & 0 & 0\\
   & & \ddots & & \\
    0 & 0 & & -b_{k+m}\iota^{(1)}_{m1} & -b_{k+m}\iota^{(1)}_{m2}
\end{array}\right]_{m \times 2m},
\end{equation}

\begin{equation}\label{m51}
M_{5,1}=\begin{bmatrix}
-\rho^{(2)}_{1,1} & -\rho^{(2)}_{1,2} & -\rho^{(2)}_{1,3} & \cdots & 0 & 0 & 0\\
 & &  & \ddots &  & \\
0 & 0 & 0 & \cdots & -\rho^{(2)}_{k1} & \rho^{(2)}_{k2} &\rho^{(2)}_{k3} \\
\end{bmatrix}_{k \times 3k},
\end{equation}

\begin{equation}\label{m52}
M_{5,2}=\begin{bmatrix}
c_{1,1} \iota^{(2)}_{1,1} & c_{1,1} \iota^{(2)}_{1,2} & \cdots & c_{1m} \iota^{(2)}_{m,1} & c_{1m} \iota^{(2)}_{m,2}\\
\vdots & \vdots & \ddots & \vdots & \vdots\\
c_{k1} \iota^{(2)}_{1,1} & c_{k1} \iota^{(2)}_{1,2} & \cdots & c_{km} \iota^{(2)}_{m,1}& c_{km} \iota^{(2)}_{m,2}
\end{bmatrix}_{k \times m},
\end{equation}

\begin{equation}\label{m61}
M_{6,1}=\begin{bmatrix}
b_{1}^{-1} \rho^{(3)}_{1,1} & b_{1}^{-1} \rho^{(3)}_{1,2} & b_{1}^{-1} \rho^{(3)}_{1,3} & \cdots & b_{k}^{-1} \rho^{(3)}_{k 1} & b_{k}^{-1} \rho^{(3)}_{k 2} & b_{k}^{-1} \rho^{(3)}_{k 3}
\end{bmatrix}_{1 \times 2k},
\end{equation}
\begin{equation}
M_{6,2}=\begin{bmatrix}\label{m62}
-b_{k+1}^{-1} \iota^{(3)}_{1,1} & -b_{k+1}^{-1} \iota^{(3)}_{1,2} & \hdots & -b_{k+m}^{-1} \iota^{(3)}_{m,1} & -b_{k+m}^{-1} \iota^{(3)}_{m,2}
\end{bmatrix}_{1 \times m},
\end{equation}

\begin{equation}\label{m71}
M_{7,1}=\begin{bmatrix}
a_{1}^{-1} \rho^{(4)}_{1,1} & a_{1}^{-1} \rho^{(4)}_{1,2} & a_{1}^{-1} \rho^{(4)}_{1,3} & \cdots & a_{k}^{-1} \rho^{(4)}_{k 1} & a_{k}^{-1} \rho^{(4)}_{k 2} & a_{k}^{-1} \rho^{(4)}_{k 3}
\end{bmatrix}_{1 \times 2k},
\end{equation}
and
\begin{equation}
M_{7,2}=\begin{bmatrix}\label{m72}
-a_{k+1}^{-1} \iota^{(4)}_{1,1} & -a_{k+1}^{-1} \iota^{(4)}_{1,2} & \hdots & -a_{k+m}^{-1} \iota^{(4)}_{m,1} & -a_{k+m}^{-1} \iota^{(4)}_{m,2}
\end{bmatrix}_{1 \times m}.
\end{equation}

All matrices, blocks, and terms $\rho_{il}^{(\vartheta)}, \iota_{jp}^{(\vartheta)}$ above are functions of parameters $\lambda_{il},\beta_{jp}$; Table \ref{termosmatriz} shows the mathematical expressions of this dependence. The numbers $a_i,a_{k+j},b_i,b_{k+j}$ and $c_{ij}$ come from the boundary conditions, with $c_{ij}$ being the generic term of the matrix $C$.

\begin{table}[H]
\begin{tabular}{ccc}
 term & & expression on parameters\\
 $\rho^{(\vartheta)}_{il}$ & & $2\sin{(\frac\pi{10})}\left[2\cos\left(-\frac{2\pi}{5}\lambda _{il}+\frac{2\pi}{5}+\frac{2\pi\vartheta}{5}\right)+1\right]$\\[4pt]

$\iota^{(\vartheta)}_{jp}$ & & $4\sin{(\frac\pi{10})}\cos\left(-\frac{\pi}{5}\beta_{jp}+\frac{\pi}{5}+\frac{\pi\vartheta}{5}\right)$ 
\end{tabular}
\caption{Expressions of matrices terms}
\label{termosmatriz}
\end{table}
Now, we can enunciate the principal result of this work, which gives a criterium to obtain a solution for the Cauchy problem \eqref{Kawahara} with boundary conditions $\eqref{condfront}$.
\begin{theorem}\label{grandeteorema}
    Assume $s\in \left(-\frac{1}{2},\frac{5}{2}\right)\setminus\left\{\frac{1}{2},\frac{3}{2}\right\}$ to be a fixed real number and $u_{0,i}\in H^s(-\infty,0)$ and $v_{0,j}\in H^s(0,+\infty)$. 
    Suppose that there exist $\lambda_{il}(s)$ and $\beta_{ps}(s)$ that satisfy 
        \begin{equation}\label{hyp}
            \displaystyle max\{s-2,0\} < \lambda_{il}(s),\beta_{ps}(s) < min\left\{s+\frac{1}{2},\frac{1}{2}\right\}
        \end{equation}
        and such that the matrix $\mathbf{M}(\boldsymbol\lambda, \boldsymbol\beta)$  is invertible. Then, there exist a time $T>0$ and a solution $(u_1,\hdots,u_k,v_1,\hdots,v_m)$ in the space $C([0,T];H^s(\mathcal{G}))$  of the Cauchy problem $\eqref{Kawahara}$ with boundary conditions $\eqref{condfront}$ satisfying the compatibility conditions \eqref{compatibility1} when $s>\frac{1}{2}$ and \eqref{compatibility2} when $s>\frac{3}{2}$. Furthermore, the data-to-solution map $(u_{0,1},\hdots,u_{0k},v_{0,1},\hdots,v_{0m})\mapsto (u_1,\hdots,u_k,v_1,\hdots,v_m)$ is locally Lipschitz continuous from $H^s(\mathcal{G})$ to $C([0, T];H^s(\mathcal{G}))$.
\end{theorem}
The principal ingredients to prove Theorem \ref{grandeteorema} are the use of boundary forcing operators, associated with the linearized Kawahara equation, constructed in \cite{Cavalcantekwak2020}, and also the ideas of the works \cite{fourth} and \cite{MC}. The main difficulty in the proof is the construction of an extended integral formula based on a Riemann-Liouville fractional interaction and involving convenient choices of some forcing functions, depending strongly on the boundary conditions, that appear in the extended problem. 

\section{Notations and function spaces} For $\phi=\phi(x)\in S(\mathbb{R})$,  $\displaystyle \hat{\phi}(\xi)=\int e^{-i\xi x}\phi(x)dx$ denotes the Fourier transform of $\phi$. For $u=u(x,t)\in S(\mathbb{R}^2)$, $$\hat{u}=\hat{u}(\xi,\tau)=\displaystyle \int e^{-i(\xi x+\tau t)}u(x,t)dxdt$$ denotes its space-time Fourier transform. We will denote  the space and the time Fourier transform of $u(x,t)$, respectively by $\mathcal{F}_{x}u(\xi,t)$ and $\mathcal{F}_{t}u(x,\tau)$.  

For any real number $\xi$ we put $\langle \xi \rangle:=1+|\xi|$ and $f(\xi,\tau) \lesssim g(\xi,\tau)$ means that there is a constant $C$ such that 
$f(\xi,\tau) \leq Cg(\xi,\tau)$ for all $(\xi, \tau)\in \R^2$. 
The characteristic function of an arbitrary set $A$ is denoted by $\chi_{A}$. Throughout the paper, we fix a cutoff function $\psi \in C_0^{\infty}(\mathbb{R})$ such that
\begin{equation}\label{cutt}
\psi(t)=
\begin{cases}
1 & \text{if}\; |t|\le 1,\\
0 & \text{if}\; |t|\ge 2
\end{cases}
\end{equation}
and $\mathbb{R}^*=\mathbb{R}\setminus \{0\}.$

\subsection{Function Spaces}

We first need to review some function spaces relevant to both the half-line and the entire real line.

Initially, for $s\geq 0$ we say that $\phi \in H^s(\mathbb{R}^+)$ if there exists $\tilde{\phi}\in H^s(\mathbb{R})$ such that
$\phi=\tilde{\phi}|_{\R+}$.  In this case we set $\|\phi\|_{H^s(\mathbb{R}^+)}:=\inf\limits_{\tilde{\phi}}\|\tilde{\phi}\|_{H^{s}(\mathbb{R})}$. For $s\geq 0$ define $$H_0^s(\mathbb{R}^+)=\Big\{\phi \in H^{s}(\mathbb{R});\,\text{supp} (\phi) \subset[0,+\infty) \Big\}.$$ For $s<0$, define $H^s(\mathbb{R}^+)$ and $H_0^s(\mathbb{R}^+)$  as the dual space of $H_0^{-s}(\mathbb{R}^+)$ and  $H^{-s}(\mathbb{R}^+)$, respectively.

Also, we define 
$$C_0^{\infty}(\mathbb{R}^+)=\Big\{\phi\in C^{\infty}(\mathbb{R});\, \text{supp}(\phi) \subset [0,+\infty)\Big\}$$
and  as those members of $C_0^{\infty}(\mathbb{R}^+)$ with compact support. We recall that $C_{0}^{\infty}(\mathbb{R}^+)$ is dense in $H_0^s(\mathbb{R}^+)$ for all $s\in \mathbb{R}$. A definition for $H^s(\R^-)$ and $H_0^s(\R^-)$ can be given analogous to that for $H^s(\R^+)$ and $H_0^s(\R^+)$.

The following results summarize useful properties of the Sobolev spaces on the half-line and their proofs could be found in \cite{CK}.

\begin{lemma}\label{sobolevh0}
	For all $f\in H^s(\mathbb{R})$  with $-\frac{1}{2}<s<\frac{1}{2}$ we have
	\begin{equation*}
	\|\chi_{(0,+\infty)}f\|_{H^s(\mathbb{R})}\lesssim  \|f\|_{H^s(\mathbb{R})}.
	\end{equation*}
\end{lemma}
\begin{lemma}\label{alta}
	If $\frac{1}{2}<s<\frac{3}{2}$ the following statements are valid:
	\begin{enumerate}
		\item [(a)] $H_0^s(\R^+)=\big\{f\in H^s(\R^+);f(0)=0\big\},$\medskip
		\item [(b)] If  $f\in H^s(\R^+)$ with $f(0)=0$, then $\|\chi_{(0,+\infty)}f\|_{H_0^s(\R^+)}\lesssim \|f\|_{H^s(\R^+)}$.
	\end{enumerate}
\end{lemma}

\begin{lemma}\label{cut}
	If $f\in  H_0^s(\mathbb{R}^+)$ with $s\in \R$, we then have
		\begin{equation*}
	\|\psi f\|_{H_0^s(\mathbb{R}^+)}\lesssim \|f\|_{H_0^s(\mathbb{R}^+)}.
	\end{equation*}
	
\end{lemma}
\begin{remark}
	In Lemmas \ref{sobolevh0}, \ref{alta} and \ref{cut}, all constants $c$ depend only on $s$ and $\psi$.
\end{remark}

Now, we introduce adpated Bourgain spaces used on \cite{Cavalcantekwak2020} in the context of half-lines.

We denote by $X^{s,b}$ the so-called Bourgain spaces associated with the linear Kawahara equation. More precisely, $X^{s,b}$  is the completion of $S'(\mathbb{R}^2)$ with respect to the norm
\begin{equation*}\label{Bourgain-norm}
\norm{f}_{X^{s,b}}^2 = \int_{\R^2} \bra{\xi}^{2s}\bra{\tau - \xi^5}^{2b}|\wt{f}(\tau,\xi)|^2 \; d\xi d\tau,
\end{equation*}
To obtain our results we also need to define the following auxiliary modified Bougain spaces introduced in  \cite{Cavalcantekwak2020}. Let  $Y^{s,b}$ and $D^{\alpha}$ be the completion of $S'(\R^2)$ with respect to the norms:
\begin{align*}
&\norm{f}_{Y^{s,b}}^2 = \int_{\R^2} \bra{\tau}^{\frac{2s}{5}}\bra{\tau - \xi^5}^{2b}|\wt{f}(\tau,\xi)|^2 \; d\xi d\tau\\\intertext{and}
&\norm{f}_{D^{\alpha}}^2 = \int_{\R^2} \bra{\tau}^{2\alpha}\mathbf{1}_{\set{\xi : |\xi| \le 1}}(\xi)|\wt{f}(\tau,\xi)|^2 \; d\xi d\tau.
\end{align*}

If $\Omega$ is a domain in $\mathbb{R}^2$, then define by $X^{s, b}(\Omega)$ and $Y^{s, b}(\Omega)$ restrictions of $X^{s, b, \beta, \sigma}$ and $Y^{s, b, \beta, \sigma}$ on $\Omega$, respectively, with natural restriction norms.

On the context of a star graph $\mathcal{G}$ we consider the Bourgain spaces $X^{s, b, \beta, \sigma}(\mathcal{G} \times (0, T))$ by

\begin{equation}
\mathcal{X}^{s, b}(\mathcal{G} \times (0, T)) 
= \bigoplus_{e \in E_-} X^{s, b}((-\infty, 0) \times (0, T)) \oplus \bigoplus_{e \in E_+} X^{s, b}((0, \infty) \times (0, T)). \tag{2.4}
\end{equation}

In a similar way, we can define the spaces $\mathcal U^{s, b}(\mathcal{G} \times (0, T))$
 and $\mathcal D^{\alpha}(\mathcal{G} \times (0, T))$.

 For solving the Cauchy problem \eqref{Kawahara}-\eqref{condfront} we will consider the following functional space $\mathcal Z^{s, b,  \alpha}(T)$, $T > 0$, $s \ge 1$,

\begin{equation}\label{functional10}
    \mathcal Z^{s, b, \alpha}(T) = \left\{ w : \mathcal{G} \times [0, T] \to \mathbb{R} : 
    \begin{aligned}
        &w \in C([0, T]; H^s(\mathcal{G})) \\
        &\cap C(\mathcal{G}; H^{\frac{s+2}{5}}([0, T])) \cap \mathcal X^{s, b}(\mathcal G\times (0,T))\cap \mathcal V^{\alpha}(\mathcal G\times (0,T)); \\
        &\partial_x^\vartheta w \in C(\mathcal{G}; H^{\frac{s+2-\vartheta}{5}}([0, T])),\; \vartheta=0,1,2,3,4
    \end{aligned}
    \right\}.
\end{equation}
 
with the following norm
\begin{equation*}
\|(u_1,\hdots,u_k,v_1,\hdots,v_m)\|_{\mathcal Z^{(s,b,\alpha)}}:=\displaystyle\sum_{i=1}^k\|u_i\|_{\mathcal Z_i^{(s,b,\alpha)}}+\displaystyle\sum_{j=1}^m\|v_j\|_{\mathcal Z_{k+j}^{(s,b,\alpha)}},
\end{equation*}
where
\begin{equation}
\begin{split}
\|w\|_{\mathcal Z_i^{(s,b,\alpha)}}=&\|w\|_{C(\R_t;H^s(\R_x))}+\sum_{\vartheta
=0}^{4}\|\partial_x^\vartheta
w\|_{C(\R_x;H^{\frac{s+2-\vartheta
}{5}}(\R_t))}+\|w\|_{ X^{s,b}}+\|w\|_{D^{\alpha}}.
\end{split}
\end{equation}

The next nonlinear estimates, in the context of the Kawahara equation for $b<\frac{1}{2}$, were derived by  \cite{Cavalcantekwak2020}  .
\begin{lemma}\label{bilinear1}
	\begin{itemize}
		\item [(a)]
		For $-7/4 < s $, there exists $b = b(s) < 1/2$ such that for all $\alpha > 1/2$, we have \begin{equation}\label{eq:bilinear1} \norm{\partial_x(uv)}_{X^{s,-b}} \le c\norm{u}_{X^{s,b} \cap D^{\alpha}}\norm{v}_{X^{s,b} \cap D^{\alpha}}. \end{equation}
		\item[(b)] For $-7/4 < s < 5/2$, there exists $b = b(s) < 1/2$ such that for all $\alpha > 1/2$, we have \begin{equation}\label{eq:bilinear2} \norm{\partial_x(uv)}_{Y^{s,-b}} \le c\norm{u}_{X^{s,b} \cap D^{\alpha}}\norm{v}_{X^{s,b} \cap D^{\alpha}}. \end{equation}
	\end{itemize}
\end{lemma}

\subsection{Riemann-Liouville fractional integral}
For Re $\alpha>0$, the tempered distribution $\frac{t_+^{\alpha-1}}{\Gamma(\alpha)}$ is defined as a locally integrable function by 
\begin{equation*}
\left \langle \frac{t_+^{\alpha-1}}{\Gamma(\alpha)},\ f \right \rangle:=\frac{1}{\Gamma(\alpha)}\int_0^{+\infty} t^{\alpha-1}f(t)dt.
\end{equation*}

For Re $\alpha>0$, integration by parts implies that
\begin{equation*}
\frac{t_+^{\alpha-1}}{\Gamma(\alpha)}=\partial_t^k\left( \frac{t_+^{\alpha+k-1}}{\Gamma(\alpha+k)}\right)
\end{equation*}
for all $k\in\mathbb{N}$. This expression allows to extend the definition, in the sense of distributions,  of $\frac{t_+^{\alpha-1}}{\Gamma(\alpha)}$ to all $\alpha \in \mathbb{C}$.

If $f\in C_0^{\infty}(\mathbb{R}^+)$, we define
\begin{equation*}
\mathcal{I}_{\alpha}f=\frac{t_+^{\alpha-1}}{\Gamma(\alpha)}*f.
\end{equation*}
Thus, for Re $\alpha>0$,
\begin{equation*}
\mathcal{I}_{\alpha}f(t)=\frac{1}{\Gamma(\alpha)}\int_0^t(t-s)^{\alpha-1}f(s)ds
\end{equation*} 
and notice that 
$$\mathcal{I}_0f=f,\quad  \mathcal{I}_1f(t)=\int_0^tf(s)ds,\quad \mathcal{I}_{-1}f=f'\quad  \text{and}\quad  \mathcal{I}_{\alpha}\mathcal{I}_{\beta}=\mathcal{I}_{\alpha+\beta}.$$

The following result states important properties of the Riemann-Liouville fractional integral operators and their proof can be found in \cite{Holmer}.

\begin{lemma}\label{lio-lemaint}
	If $0\leq \alpha <\infty$,\, $s\in \mathbb{R}$ and $\varphi \in C_0^{\infty}(\mathbb{R})$, then we have
	\begin{align}
	&\|\mathcal{I}_{-\alpha}h\|_{H_0^s(\mathbb{R}^+)}\leq c \|h\|_{H_0^{s+\alpha}(\mathbb{R}^+)}\label{lio}\\
	\intertext{and}
	&\|\varphi \mathcal{I}_{\alpha}h\|_{H_0^s(\mathbb{R}^+)}\leq c_{\varphi} \|h\|_{H_0^{s-\alpha}(\mathbb{R}^+)}.\label{lemaint}
	\end{align}
\end{lemma}
\section{Estimates associated of some operators}\label{section3}

\subsection{Free popagator} The  unitary  group associated to the linear Kawahara equation is defined as
\begin{equation*}
e^{t\partial_x^5}\phi(x)=\frac{1}{2\pi}\int e^{ix\xi}e^{it\xi^5}\hat{\phi}(\xi)d\xi,
\end{equation*}
that satisfies
\begin{equation}\label{lineark}
\begin{cases}
(\partial_t-\partial_x^5)e^{t\partial_x^5}\phi(x) =0,& (t,x)\in\mathbb{R}\times\mathbb{R},\\
e^{t\partial_x^5}\phi(x)\big|_{t=0}=\phi(x),& x\in\mathbb{R}.
\end{cases}
\end{equation}

We will establish a series of estimates for the solutions of linear IVP \eqref{lineark}. We refer the reader to \cite{Cavalcantekwak2020} for a proof of these results.

\begin{lemma}\label{grupok}
	Let $s\in\mathbb{R}$ and $b, \alpha \in \R$. If $\phi\in H^s(\mathbb{R})$, then
	\begin{enumerate}
		\item[(a)] \emph{(Space traces)} $\|\psi(t)e^{t\partial_x^5}\phi(x)\|_{C_t\left(\mathbb{R};\,H_x^s(\mathbb{R})\right)}\lesssim \|\phi\|_{H^s(\mathbb{R})}$;
		\item[(b)] \emph{((Derivatives) Time traces)} 
		\[
		\|\psi(t) \partial_x^{j}e^{t\partial_x^5}\phi(x)\|_{C_x\left(\mathbb{R};H_t^{\frac{s+2-j}{5}}(\mathbb{R})\right)}\lesssim_{\psi, s, j} \|\phi\|_{H^s(\mathbb{R})}, \quad j\in\{0,1,2,3,4\};
		\]
		\item [(c)] \emph{(Bourgain spaces)} $\|\psi(t)e^{t\partial_x^5}\phi(x)\|_{X^{s,b}\cap D^{\alpha}}\lesssim_{\psi, b, \alpha}  \|\phi\|_{H^s(\mathbb{R})}$.
	\end{enumerate}
\end{lemma}
\begin{remark}
	The spaces $D^{\alpha}$ introduced in \cite{Holmer}  give us  useful auxiliary norms of the classical Bourgain spaces in order to validate the bilinear estimates associated to the Kawahara equation for $b<\frac12$ and $\alpha>\frac12$ (see Lemma \ref{bilinear1}).  
\end{remark}

\subsection{The Duhamel boundary forcing operator associated to the linear Kawahara equation}\label{section4}
Now we give the properties of the Duhamel boundary forcing operator introduced in \cite{Cavalcantekwak2020}, that is
\begin{equation}\label{vk}
\begin{split}
\mathcal{L}^0f(t,x)&=M\int_0^te^{(t-t')\partial_x^5}\delta_0(x)\mathcal{I}_{-\frac45}f(t')dt'\\
	&=M\int_0^t B\left(\frac{x}{(t-t')^{1/5}}\right)\frac{\mathcal{I}_{-\frac45}f(t')}{(t-t')^{1/5}}dt',
\end{split}
\end{equation}
defined for all $f\in C_0^{\infty}(\mathbb{R}^+)$ with \begin{equation}\label{eq:M}
M = \dfrac{1}{B(0)\Gamma(4/5)},
\end{equation} and $B$ denotes the Oscillatory integral 
$$B(x)=\frac{1}{2\pi}\int_{\R}e^{ix\xi}e^{i \xi^5} \; d\xi.$$
From definition of $\mathcal{L}^0$, it follows 
\begin{equation}\label{forcingl}
\begin{cases}
(\partial_t-\partial_x^5)\mathcal{L}^0f(t,x)=M\delta_0(x)\mathcal{I}_{-\frac{4}{5}}f(t)& \text{for}\quad (x,t)\in \mathbb{R}\times\mathbb{R},\\
\mathcal{L}^0f(0,x)=0& \text{for}\quad x\in\mathbb{R}.
\end{cases}
\end{equation}

The proof of the results exhibited in this section was shown in \cite{Cavalcantekwak2020}.

\begin{lemma}[Continuity of $\partial_x^k \mathcal{L}^0f(t,x)$]\label{lemacrb1}
	Let $f\in C_{0}^{\infty}(\mathbb{R}^+)$.
	\begin{itemize}
\item[(a)] For fixed $ 0 \le t \le 1$, $\partial_x^k \mathcal{L}^0f(t,x)$, $k=0,1,2,3$, is continuous in $x \in \mathbb{R}$ and has the decay property in terms of the spatial variable as follows:
\begin{equation}\label{eq:decay1}
|\partial_x^k \mathcal{L}^0f(t,x)| \lesssim_{N} \norm{f}_{H^{N+k}}\bra{x}^{-N}, \quad N \ge 0.
\end{equation}
\item[(b)] For fixed $ 0 \le t \le 1$, $\partial_x^4\mathcal{L}^0f(t,x)$ is continuous in $x$ for $x\neq 0$ and is discontinuous at $x =0$ satisfying
\[\lim_{x\rightarrow 0^{-}}\partial_x^4\mathcal{L}^0f(t,x)=c_1\mathcal{I}_{-4/5}f(t),\ \lim_{x\rightarrow 0^{+}}\partial_x^4\mathcal{L}^0f(t,x)=c_2\mathcal{I}_{-4/5}f(t)\]
for $c_1 \neq c_2$. $\partial_x^4\mathcal{L}^0f(t,x)$ also has the decay property in terms of the spatial variable
\begin{equation}\label{eq:decay2}
|\partial_x^4\mathcal{L}^0f(t,x)| \lesssim_{N} \norm{f}_{H^{N+4}}\bra{x}^{-N}, \quad N \ge 0.\\
\end{equation}
	\end{itemize}
\end{lemma}
\medskip

From \eqref{eq:M} and  \eqref{vk}, we have that $\mathcal{L}^0f(t,0)=f(t).$
Also, it follows from \eqref{vk} (for details see \cite{Cavalcantekwak2020}) the function 
\begin{equation}
v(x,t)=e^{t\partial_x^5}\phi(x)+\mathcal{V}\big(f-e^{t\partial_x^5}\phi\big|_{x=0}\big)(x,t),
\end{equation}
where $f\in C_0^{\infty}(\R^+)$ and $\phi \in S(\R)$, solves the linear problem
\begin{equation}\label{linearr}
\begin{cases}
(\partial_t-\partial_x^5)v(x,t)=0& \text{for}\quad (x,t)\in \R^*\times\mathbb{R},\\
v(x,0)=\phi(x)& \text{for}\quad x\in\mathbb{R},\\
v(0,t)=f(t)& \text{for}\quad t\in(0, +\infty),
\end{cases}
\end{equation}
in the sense of distributions.

\subsection{The Duhamel Boundary Forcing Operator Classes associated to linear KdV equation}

To extend our results to a broader range of regularity, we utilize two classes of boundary forcing operators, generalizing the operator $\mathcal{L}^0$ introduced in \cite{Cavalcantekwak2020}.

Let $\lambda\in \mathbb{C}$ with 
$\text{Re}\,\lambda>-3$ and $g\in C_0^{\infty}(\mathbb{R}^+)$. Define the operators
\begin{equation*}
\mathcal{L}_{\pm}^{\lambda}g(t,x)=\left[\frac{x_{\mp}^{\lambda-1}}{\Gamma(\lambda)}*\mathcal{L}^0\big(\mathcal{I}_{-\frac{\lambda}{5}}g\big)(t, \cdot)   \right](x),
\end{equation*}
with $\frac{x_{-}^{\lambda-1}}{\Gamma(\lambda)}=e^{i\pi \lambda}\frac{(-x)_{+}^{\lambda-1}}{\Gamma(\lambda)}$. Then, 
using \eqref{forcingl} we have that
\begin{equation*}
(\partial_t-\partial_x^5)\mathcal{L}_{-}^{\lambda}g(t,x)=M\frac{x_{+}^{\lambda-1}}{\Gamma(\lambda)}\mathcal{I}_{-\frac{4}{5}-\frac{\lambda}{5}}g(t)
\end{equation*}
and
\begin{equation*}
(\partial_t-\partial_x^5)\mathcal{L}_{+}^{\lambda}g(t,x)=M\frac{x_{-}^{\lambda-1}}{\Gamma(\lambda)}\mathcal{I}_{-\frac{4}{5}-\frac{\lambda}{5}}g(t).
\end{equation*}

The following lemmas state properties of the operators classes $\mathcal{L}_{\pm}^{\lambda}$. For the proofs
we refer the reader \cite{Cavalcantekwak2020}.
\begin{lemma}[Spatial continuity and decay properties for $\mathcal{L}_{\pm}^{\lambda}g(t,x)$]\label{holmer1}
	Let $g\in C_0^{\infty}(\mathbb{R}^+)$ and $M$ be as in \eqref{eq:M}. Then, we have
	\begin{equation}\label{eq:relation}
	\mathcal{L}_{\pm}^{-k}g=\partial_x^k\mathcal{L}^{0}\mathcal{I}_{\frac{k}{5}}g, \qquad k=0,1,2,3,4.
	\end{equation}
	Moreover, $\mathcal{L}_{\pm}^{-4}g(t,x)$ is continuous in $x \in \R \setminus \set{0}$ and has a step discontinuity of size $Mg(t)$ at $x=0$. 
    \end{lemma}

\begin{lemma}[Values of $\mathcal{L}_{+}^{\lambda}f(t,0)$ and $\mathcal{L}_{-}^{\lambda}f(t,0)$ ]\label{tracel}
	 For $\mbox{Re}\ \lambda>-4$,
	\begin{equation}\label{lr0}
	\mathcal{L}_{+}^{\lambda}f(t,0)=\frac{1}{B(0)\Gamma(4/5)}\frac{\cos\left(\frac{(1+4\lambda)\pi}{10}\right)}{5\sin\left(\frac{(1-\lambda)\pi}{5}\right)}f(t)=4\sin\left(\frac{\pi}{10}\right)\cos\left(\dfrac{\pi(1-\lambda)}{5}\right)f(t)
	\end{equation}
	and
	\begin{equation}\label{ll0}
\mathcal{L}_{-}^{\lambda}f(t,0)= \frac{1}{B(0)\Gamma(4/5)}\frac{\cos\left(\frac{(1-6\lambda)\pi}{10}\right)}{5\sin\left(\frac{(1-\lambda)\pi}{5}\right)}f(t)=2\sin\left(\frac{\pi}{10}\right)\left[2\cos\left(\dfrac{2\pi(1-\lambda)}{5}\right)+1\right]f(t)
	\end{equation}
\end{lemma}
\begin{lemma}\label{edbf}$\;$\\
	\begin{itemize}
\item[(a)] \emph{(Space traces)}  Let $-\frac92 < s < 5$. For $\max(s-\frac{9}{2}, -4) <\lambda< \min(s+\frac{1}{2}, \frac12)$, we have
		$$\|\psi(t)\mathcal{L}_{\pm}^{\lambda}f(t,x)\|_{C\big(\mathbb{R}_t;\,H^s(\mathbb{R}_x)\big)}\leq c \|f\|_{H_0^\frac{s+2}{5}(\mathbb{R}^+)};$$
\item[(b)] \emph{((Derivatives) Time traces)} For $-4+j<\lambda<1+j$, $j=0,1,2,3,4$, we have 
		\begin{equation}\label{eq:(b)0}
		\|\psi(t)\partial_x^j\mathcal{L}_{\pm}^{\lambda}f(t,x)\|_{C\big(\mathbb{R}_x;\,H_0^{\frac{s+2-j}{5}}(\mathbb{R}_t^+)\big)}\leq c \|f\|_{H_0^\frac{s+2}{5}(\mathbb{R}^+)};
		\end{equation}
\item[(c)] \emph{(Bourgain spaces)} Let $-7 < s < \frac52$ and $b < \frac12 < \alpha < 1-b$. For $\max(s-2, -\frac{13}{2}) <\lambda< \min(s+\frac{1}{2}, \frac12)$, we have
		\[\|\psi(t)\mathcal{L}_{\pm}^{\lambda}f(t,x)\|_{X^{s,b}\cap D^{\alpha}}\leq c \|f\|_{H_0^\frac{s+2}{5}(\mathbb{R}^+)}.\]
\end{itemize}
\end{lemma}
\begin{remark}\label{remark} The range of regularity where our estimates hold varies depending on the index $\lambda$. For example, when $\lambda = 0$, the Bourgain spaces estimates are valid only for  $-1/2 < s < 2$. This dependence on $\lambda$ is one of the reasons why we must consider the entire family of operators, rather than focusing solely on $\mathcal{L}^0$.
\end{remark}
\subsection{{The Duhamel Inhomogeneous Solution Operator}}\label{sectionduhamel}

We introduce the Duhamel inhomogeneous solution operator $\mathcal{D}$ as
\begin{equation*}
\mathcal{D}w(t,x)=\int_0^te^{(t-t')\partial_x^5}w(t',x)dt',
\end{equation*}
it follows that
\[\left \{
\begin{array}{l}
(\partial_t-\partial_x^5)\mathcal{D}w(t,x) =w(t,x),\ (t,x)\in\mathbb{R}\times\mathbb{R},\\
\mathcal{D}w(x,0) =0,\ x\in\mathbb{R}.
\end{array}
\right.\]

The following lemma summarizes useful estimates for the Duhamel inhomogeneous solution operators  $\mathcal{D}$  that will be used later in the proof of the main results and its proof can be seen in \cite{Cavalcantekwak2020}.

\begin{lemma}\label{duhamelk} For all $s\in\mathbb{R}$ we have the following estimates: \medskip 
		Let $s\in \mathbb{R}$. For $0 < b < \frac12 < \alpha < 1-b$, we have
	\begin{itemize}
		\item[(a)] \emph{(Space traces)} 
		\begin{equation*}
		\|\psi(t)\mathcal{D}w(t,x)\|_{C\big(\mathbb{R}_t;\,H^s(\mathbb{R}_x)\big)}\lesssim \|w\|_{X^{s,-b}};
		\end{equation*}
		\item[(b)] \emph{((Derivatives) Time traces)}
		\begin{equation*}
		\|\psi(t) \partial_x^j\mathcal{D}w(t,x)\|_{C(\mathbb{R}_x;H^{\frac{s+2-j}{5}}(\mathbb{R}_t))}\lesssim	\|w\|_{X^{s,-b}}+\|w\|_{Y^{s,-b}};
		\end{equation*}
		\item[(c)] \emph{(Bourgain spaces estimates)} 
		\begin{equation*}
		\|\psi(t)\mathcal{D}w(t,x)\|_{X^{s,b} \cap D^{\alpha}}\lesssim  \|w\|_{X^{s,-b}}.
        \end{equation*}
    \end{itemize}
    \end{lemma}
Remark in (b) that $\|\psi(t) \partial_x^j\mathcal{D}w(t,x)\|_{C(\mathbb{R}_x;H_0^{\frac{s+2-j}{5}}(\mathbb{R}_t^+))}$ has same bound for $s < \frac{11}{2} + j$. 
\medskip
\section*{Notation}
In the following sections, the indices $i,j,\vartheta,l$ and $p$ are integers. $i$ is between $1$ and $k$; $j$ is between $1$ and $m$; $\vartheta$ is between $0$ and $4$; $l$ is between $1$ and $3$; and $p$ is $1$ or $2$.
\section{Start of the Proof of  Theorem \ref{grandeteorema}}\label{startproof}
Here, we present the first part of the proof of the main result.

Firstly, we will obtain an integral equation that solves the Cauchy problem \eqref{Kawahara}.
We start rewriting the first two vertex condition \eqref{condfront}  in terms of matrices as:
\begin{equation}\label{m1}
\left[\begin{array}{cccccccc}
1&-a_2& 0&\cdots&0&0&\cdots&0\\
1& 0&-a_3&\cdots&0&0&\cdots&0\\
\vdots&\vdots&\vdots&\ddots&\vdots&\vdots&\ddots&\vdots\\
1&0&0&\cdots&-a_k&0&\cdots&0\\
1&0&0&\cdots&0&-a_{k+1}&\cdots&0\\
\vdots&\vdots&\vdots&\ddots&\vdots&\vdots&\ddots&\vdots\\
1&0&0&\cdots&0&0&\cdots&-a_{k+m}
\end{array}\right]_{(k+m-1) \times (k+m)}
\left[\begin{array}{c}
u_1(0,t)\\
\vdots\\
u_k(0,t)\\
v_1(0,t)\\
\vdots\\
v_m(0,t)
\end{array}\right]_{(k+m)\times 1}=0
\end{equation}
and
\begin{equation}\label{m2}
\left[\begin{array}{cccccccc}
1&-b_2& 0&\cdots&0&0&\cdots&0\\
1& 0&-b_3&\cdots&0&0&\cdots&0\\
\vdots&\vdots&\vdots&\ddots&\vdots&\vdots&\ddots&\vdots\\
1&0&0&\cdots&-b_k&0&\cdots&0\\
1&0&0&\cdots&0&-b_{k+1}&\cdots&0\\
\vdots&\vdots&\vdots&\ddots&\vdots&\vdots&\ddots&\vdots\\
1&0&0&\cdots&0&0&\cdots&-b_{k+m}
\end{array}\right]_{(k+m-1) \times (k+m)}
\left[\begin{array}{c}
\partial_x u_1(0,t)\\
\vdots\\
\partial_x u_k(0,t)\\
\partial_x v_1(0,t)\\
\vdots\\
\partial_x v_m(0,t)
\end{array}\right]_{(k+m)\times 1}=0.
\end{equation}.

Now, if we write $C$ as 
\begin{equation*}\label{C}
\left[\begin{array}{cccc}
c_{1,1}& c_{1,2}& \cdots&c_{1m}\\
c_{2,1}& c_{2,2}& \cdots&c_{2m}\\
\vdots& \vdots& \ddots&\vdots\\
c_{k1}& c_{k2}& \cdots&c_{km}\\
\end{array}\right]_{k\times m},
\end{equation*}

then the third condition on \eqref{condfront} becomes
\begin{equation}\label{m3}
\left[\begin{array}{ccccccc}
-1&0&\cdots&0&c_{1,1}&  \cdots&c_{1m}\\
0&-1&\cdots&0&c_{2,1}& \cdots&c_{2m}\\
\vdots&\vdots&\ddots&\vdots&  \vdots&\ddots&\vdots\\0&0&\cdots&-1&c_{k1}& \cdots&c_{km}\\
\end{array}\right]_{k\times (k+m)}
\left[\begin{array}{c}
\partial^2_{x} u_1(0,t)\\
\vdots\\
\partial^2_{x} u_k(0,t)\\
\partial^2_{x} v_1(0,t)\\
\vdots\\
\partial^2_{x} v_m(0,t)\\
\end{array}\right]_{(k+m)\times 1}=0.
\end{equation}

Finally, the last two conditions are
\begin{equation}\label{m4}
\left[\begin{array}{cccccc}
b_1^{-1}& \cdots& b_k^{-1}&-b_{k+1}^{-1}& \cdots& -b_{k+m}^{-1}\\
\end{array}\right]_{1\times (k+m)}
\left[\begin{array}{c}
\partial^3_{x} u_1(0,t)\\
\vdots\\
\partial^3_{x} u_k(0,t)\\
\partial^3_{x} v_1(0,t)\\
\vdots\\
\partial^3_{x} v_m(0,t)\\
\end{array}\right]_{(k+m)\times 1}=0
\end{equation}
and

\begin{equation}\label{m5}
\left[\begin{array}{cccccc}
a_1^{-1}& \cdots& a_k^{-1}&-a_{k+1}^{-1}& \cdots& -a_{k+m}^{-1}\\
\end{array}\right]_{1\times (k+m)}
\left[\begin{array}{c}
\partial^4_{x} u_1(0,t)\\
\vdots\\
\partial^4_{x} u_k(0,t)\\
\partial^4_{x} v_1(0,t)\\
\vdots\\
\partial^4_{x} v_m(0,t)\\
\end{array}\right]_{(k+m)\times 1}=0.
\end{equation}

Let $\widetilde{u}_{0i}$ and $\widetilde{v}_{0j}$ be nice extensions of ${u}_{0i}$ and ${v}_{0j}$, respectively satisfying
\begin{equation*}
\|\widetilde{u}_{0i}\|_{H^s(\R)}\leq c\|{u}_{0i}\|_{H^s(\R^{+})} \text{ and }\; \|\widetilde{v}_{0j}\|_{H^s(\R)}\leq c\|{v}_{0j}\|_{H^s(\R^{+})}.
\end{equation*}
Initially, based in the work \cite{Cavalcantekwak2020}, we look for solutions in the form
\begin{align*}
&{u}_{i}(x,t)=\mathcal{L}_-^{\lambda_{i1}}\gamma_{i1}(x,t)+\mathcal{L}_-^{\lambda_{i2}}\gamma_{i2}(x,t)+\mathcal{L}_-^{\lambda_{i3}}\gamma_{i3}(x,t)+F_i(x,t),\\
&{v}_{j}(x,t)=\mathcal{L}_+^{\beta_{j1}}\theta_{j1}(x,t)+\mathcal{L}_+^{\beta_{j2}}\theta_{j2}(x,t)+G_j(x,t),
\end{align*}
where $\gamma_{i1}, \gamma_{i2}, \gamma_{i3}, \theta_{j1}$  and $\theta_{j2}$ are unknown functions that will be choosen later, and
\begin{align*}
&F_i(x,t)=e^{t\partial_x^5}\widetilde{u}_{0i}+\mathcal{D}(\partial_x(u_i^2))(x,t),\\ 
&G_j(x,t)=e^{t\partial_x^5}\widetilde{v}_{0j}+\mathcal{D}(\partial_x(v_j^2))(x,t),
\end{align*}
where $\mathcal{D}$ is the inhomogeneous solution operator associated with the Kawahara equation.

Assuming that $\gamma_{i1},\gamma_{i2},\gamma_{i3},\theta_{j1}$ and $\theta_{j2}$ belong to $C_{0}^{\infty}(\mathbb{R}^+)$, we can apply Lemma \ref{tracel} to obtain the following trace values: 
\begin{align}
&\displaystyle \partial_x^\vartheta u_i (0,t)=2\sin\left(\frac{\pi}{10}\right)\sum_{l=1}^{3} \left[2\cos\left(-\frac{2\pi}{5}\lambda _{il}+\frac{2\pi}{5}+\frac{2\pi\vartheta}{5}\right)+1\right]\gamma_{il}(t)+F_i(0,t)\label{tec1}\quad \text{ and}\\
&\displaystyle \partial_x^\vartheta v_j (0,t)=4\sin\left(\frac{\pi}{10}\right)\sum_{p=1}^{2} \left[ \cos\left(-\frac{\pi}{5}\beta_{jp}+\frac{\pi}{5}+\frac{\pi\vartheta}{5}\right)\right]\theta_{jp}(t)+G_j(0,t).\label{tec2}
\end{align}

Using \eqref{tec1} and \eqref{tec2}, we can apply the operator $\mathcal{I}_{\frac{\vartheta}{5}}$ to rewrite \eqref{m1}-\eqref{m5} as \begin{equation}\label{mO}
M^{(\vartheta)}
\left[\begin{array}{cc}
P^{(\vartheta)} & 0\\
0& I^{(\vartheta)}\\
\end{array}\right]_{(k+m)\times(3k+2m)}\\
\left[\begin{array}{c}
\gamma_{1,1}\\
\gamma_{1,2}\\
\gamma_{1,3}\\
\vdots\\
\gamma_{k1}\\
\gamma_{k2}\\
\gamma_{k3}\\
\theta_{1,1}\\
\theta_{1,2}\\
\vdots\\
\theta_{m1}\\
\theta_{m2}\\
\end{array}\right]_{(3k+2m)\times 1}
=-M^{(\vartheta)}\left[\begin{array}{c}
\partial_x^{\vartheta} \mathcal{I}_{\frac{\vartheta}{5}} F_1(0,t)\\
\vdots\\
\partial_x^{\vartheta} \mathcal{I}_{\frac{\vartheta}{5}} F_k(0,t)\\
\partial_x^{\vartheta} \mathcal{I}_{\frac{\vartheta}{5}} G_1(0,t)\\
\vdots\\
\partial_x^{\vartheta} \mathcal{I}_{\frac{\vartheta}{5}} G_m(0,t)
\end{array}\right]_{(k+m)\times 1},
\end{equation}

where  $M^{(\vartheta)}$ are the matrices on identities \eqref{m1}-\eqref{m5} whose order is given by
\begin{equation*}
\begin{cases}
    (k+m-1)\times(k+m),& \text{ if }\ \vartheta\in \{0,1\}\ ;\\
    k\times(k+m),& \text{ if }\ \vartheta=2;\\
    1\times(k+m),& \text{ if }\ \vartheta\in\{3,4\};
\end{cases}    
\end{equation*}\\
\begin{equation}
P^{(\vartheta)}=
\left[\begin{array}{ccccccccccc}
\rho^{(\vartheta)}_{1,1}&\rho^{(\vartheta)}_{1,2}&\rho^{(\vartheta)}_{1,3}&0& 0&0& \hdots &0&0&0\\
0& 0& 0&\rho^{(\vartheta)}_{2,1}&\rho^{(\vartheta)}_{2,2}& \rho^{(\vartheta)}_{2,3}& &0&0&0\\
\vdots &\vdots &\vdots & & & & \ddots& & & \\
0&0&0&0&0& 0&&\rho^{(\vartheta)}_{k1}&\rho^{(\vartheta)}_{k2}& \rho^{(\vartheta)}_{k3}
\end{array}\right]_{(k\times 3k)};\end{equation}
and
\begin{equation}
\ I^{(\vartheta)}=
\left[\begin{array}{cccccccccccccccc}
  \iota^{(\vartheta)}_{1,1}& \iota^{(\vartheta)}_{1,2}& 0&0& \hdots&0&0\\
0& 0& \iota^{(\vartheta)}_{2,1}& \iota^{(\vartheta)}_{2,2}& &0&0\\
\vdots&\vdots&&& \ddots&&\\
0&0&0&0& & \iota^{(\vartheta)}_{m1}& \iota^{(\vartheta)}_{m2}
\end{array}\right]_{(m\times 2m)},
\end{equation}
with $\rho^\vartheta_{il}=2\sin\left(\frac{\pi}{10}\right)\left[2\cos\left(-\frac{2\pi}{5}\lambda _{il}+\frac{2\pi}{5}+\frac{2\pi\vartheta}{5}\right)+1\right]$ and $\iota^\vartheta_{jp}=4\sin\left(\frac{\pi}{10}\right)\cos\left(-\frac{\pi}{5}\beta_{jp}+\frac{\pi}{5}+\frac{\pi\vartheta}{5}\right)$.\\

Now, multiplying the matrices in \eqref{mO}, we infer
\begin{equation}\label{M}
\mathbf{M}(\boldsymbol\lambda, \boldsymbol\beta)\left[\begin{array}{c}
\boldsymbol\gamma\\
\boldsymbol\theta
\end{array}\right]
=\mathbf{F},
\end{equation}
where
\begin{equation}\label{M1}
\mathbf{M}(\boldsymbol\lambda, \boldsymbol\beta)=\left[\begin{array}{cc}
M_{1,1} & 0\\
M_{2,1} & M_{2,2}\\
\hdashline
M_{3,1} & 0\\
M_{41} & M_{42}\\
\hdashline
M_{51} & M_{52}\\
\hdashline
M_{61} & M_{62}\\
\hdashline
M_{71} & M_{72}
\end{array}\right]_{(3k+2m)\times (3k+2m)},
\end{equation}
\begin{equation}\label{M2}
\mathbf{F}=\left[\begin{array}{c}
-F_1(0,t)+a_2F_2(0,t)\\
-F_1(0,t)+a_3F_3(0,t)\\
\vdots\\
-F_1(0,t)+a_kF_k(0,t)\\
-F_1(0,t)+a_{k+1}G_1(0,t)\\
\vdots\\
-F_1(0,t)+a_{k+m}G_m(0,t)\\
\hdashline
-\partial_x \mathcal{I}_{\frac{1}{5}}F_1(0,t)+b_2\partial_x \mathcal{I}_{\frac{1}{5}}F_2(0,t)\\
-\partial_x \mathcal{I}_{\frac{1}{5}}F_1(0,t)+b_3\partial_x \mathcal{I}_{\frac{1}{5}}F_3(0,t)\\
\vdots\\
-\partial_x \mathcal{I}_{\frac{1}{5}}F_1(0,t)+b_k\partial_x \mathcal{I}_{\frac{1}{5}}F_k(0,t)\\
-\partial_x \mathcal{I}_{\frac{1}{5}}F_1(0,t)+b_{k+1}\partial_x \mathcal{I}_{\frac{1}{5}}G_1(0,t)\\
\vdots\\
-\partial_x \mathcal{I}_{\frac{1}{5}}F_1(0,t)+b_{k+m}\partial_x \mathcal{I}_{\frac{1}{5}}G_m(0,t)\\
\hdashline
\partial_x^2 \mathcal{I}_{\frac{2}{5}}F_1(0,t)-c_{1,1}\partial_x^2 \mathcal{I}_{\frac{2}{5}}G_1(0,t)-\hdots-c_{1m}\partial_x^2 \mathcal{I}_{\frac{2}{5}}G_m(0,t)\\
\partial_x^2 \mathcal{I}_{\frac{2}{5}}F_2(0,t)-c_{2,1}\partial_x^2 \mathcal{I}_{\frac{2}{5}}G_1(0,t)-\hdots-c_{2m}\partial_x^2 \mathcal{I}_{\frac{2}{5}}G_m(0,t)\\
\vdots\\
\partial_x^2 \mathcal{I}_{\frac{2}{5}}F_k(0,t)-c_{k1}\partial_x^2 \mathcal{I}_{\frac{2}{5}}G_1(0,t)-\hdots-c_{km}\partial_x^2 \mathcal{I}_{\frac{2}{5}}G_m(0,t)\\
\hdashline
-\displaystyle\sum_{i=1}^k b_i^{-1}\partial_x^3 \mathcal{I}_{\frac{3}{5}}F_i(0,t) + \sum_{j=1}^m b_{k+j}^{-1}\partial_x^3 \mathcal{I}_{\frac{3}{5}}G_j(0,t)\\
\hdashline
-\displaystyle\sum_{i=1}^k a_i^{-1}\partial_x^4 \mathcal{I}_{\frac{4}{5}}F_i(0,t) + \sum_{j=1}^m a_{k+j}^{-1}\partial_x^4 \mathcal{I}_{\frac{4}{5}}G_j(0,t)
\end{array}\right]_{(3k+2m)\times 1},
\end{equation}
and
\begin{equation}\label{M3}
\left[\begin{array}{c}
\boldsymbol\gamma\\
\boldsymbol\theta
\end{array}\right]=\left[\begin{array}{r}
\gamma_{1,1}\\
\gamma_{1,2}\\
\gamma_{1,3}\\
\vdots\\
\gamma_{k1}\\
\gamma_{k2}\\
\gamma_{k3}\\
\theta_{1,1}\\
\theta_{1,2}\\
\vdots\\
\theta_{m1}\\
\theta_{m2}\\
\end{array}\right]_{(3k+2m)\times 1}.
\end{equation}
The blocks $M_{1,1}$, $M_{2,1}$, $M_{2,2}$, $M_{3,1}$, $M_{4,1}$, $M_{4,2}$, $M_{5,1}$, $M_{5,2}$, $M_{6,1}$, $M_{6,2}$, $M_{7,1}$ and $M_{7,2}$ in \eqref{M1} are the matrices $\eqref{m11}$ - $\eqref{m72}$, respectively.\\

\begin{remark}
    The dashed lines in \eqref{M1} and \eqref{M2} divide the matrices into five blocks, whose origin is from the five equations given in \eqref{mO}.
\end{remark}

Finally, we need to obtain functions $\gamma_{il}$ and $\theta_{jp}$, and parameters $\lambda_{il}$ and $\beta_{jp}$ that satisfy \eqref{M}, since the above blocks depend on $\rho^{(\vartheta)}_{il}\text{ and } \iota_{jp}^{(\vartheta)}$, which are functions on $\lambda_{il}$'s and $\beta_{jp}$'s.\\

By using the hypothesis of Theorem \ref{grandeteorema} we fix parameters $\lambda_{il}$ and $\beta_{jp}$, such that
 \begin{equation}
\max\{s-2,0\}<\lambda_{il}(s),\beta_{jp}(s)<\min\left\{s+\frac12,\frac12\right\},
 \end{equation}
and the matrix $\mathbf{M}(\boldsymbol\lambda, \boldsymbol\beta)$ is invertible.

\section{Truncated integral operator}

In this part of the proof, we define the truncated integral operator and the appropriate function space with the goal of using the Fourier restriction method of Bourgain \cite{Bourgain1}.

Given $s$ as in the hypothesis of Theorem \ref{grandeteorema} we fix the parameters $\lambda_{il}$ and $\beta_{jp}$ as well as the functions $\gamma_{il}$ and $\theta_{jp}$ chosen as in section \ref{startproof}. Let $b=b(s)<1/2$ and $\alpha(b)>1/2$ be such that the estimates given in Lemma \ref{bilinear1} are valid.

Define the operator
\begin{equation}
\Lambda=(\Lambda^-_1,\hdots,\Lambda^-_k,\Lambda^+_1,\hdots,\Lambda^+_m)
\end{equation}
where
\begin{align*}
&\Lambda^-_i u_i(x,t)=\psi(t)\mathcal{L}_-^{\lambda_{i1}}\gamma_{i1}(x,t)+\psi(t)\mathcal{L}_-^{\lambda_{i2}}\gamma_{i2}(x,t)+\psi(t)\mathcal{L}_-^{\lambda_{i3}}\gamma_{i3}(x,t)+F_i(x,t),\\
&\Lambda^+_jv_j(x,t)=\psi(t)\mathcal{L}_+^{\beta_{j1}}\theta_{j1}(x,t)+\psi(t)\mathcal{L}_+^{\beta_{j2}}\theta_{j2}(x,t)+G_j(x,t),
\end{align*}
where
\begin{align*}
&F_i(x,t)=\psi(t)(e^{t\partial_x^5}\widetilde{u}_0+\mathcal{D}(\partial_x(u_i^2))(x,t)),\\
&G_j(x,t)=\psi(t)(e^{t\partial_x^5}\widetilde{v}_0+\mathcal{D}(\partial_x(v_j^2))(x,t)).
\end{align*}

Here, we recall the cuttoff function $\psi$ defined in subsection \ref{cutt}, which is essential to use the Fourier restriction method of Bourgain \cite{Bourgain1}.

We consider $\Lambda$ on the Banach space $\displaystyle Z_{(s,b,\alpha)}=\displaystyle\prod_{i=1}^k Z^i_{(s,b,\alpha)}\cdot\prod_{j=1}^m Z^{k+j}_{(s,b,\alpha)}$, where
\begin{equation*}
\begin{split}
Z^n_{(s,b,\alpha)}= \{w\in C(\R_t;&H^s(\R_x)) \cap C(\R_x;H^{\frac{s+2}{5}}(\R_t))\cap X^{s,b}\cap D^{\alpha};\\
&\partial_x^\vartheta
  w\in C(\R_x;H^{\frac{s+2-\vartheta
}{5}}(\R_t)), \forall\vartheta
\in\{0,1,2,3,4\}\}\; (n=1,\hdots,k+m),
\end{split}
\end{equation*}
with norm
\begin{equation*}
\|(u_1,\hdots,u_k,v_1,\hdots,v_m)\|_{Z_{(s,b,\alpha)}}=\displaystyle\sum_{i=1}^k\|u_i\|_{Z^i_{(s,b,\alpha)}}+\displaystyle\sum_{j=1}^m\|v_j\|_{Z^{k+j}_{(s,b,\alpha)}},
\end{equation*}
where
\begin{equation}
\begin{split}
\|w\|_{Z^i_{(s,b,\alpha)}}=&\|w\|_{C(\R_t;H^s(\R_x))}+\sum_{\vartheta
=0}^{4}\|\partial_x^\vartheta
w\|_{C(\R_x;H^{\frac{s+2-\vartheta
}{5}}(\R_t))}+\|w\|_{ X^{s,b}}+\|w\|_{D^{\alpha}}.\\
\end{split}
\end{equation}

\begin{remark}
 Note that the spaces $Z_{(s,b,\alpha)}$ defined previously consists of functions with domain $\prod_{n=1}^{k+m}\mathbb R \times (0,T)$, unlike the space $\mathcal Z^{s,b,\alpha}$ defined in \eqref{functional10}, where the functions have domain on $\mathcal G \times (0,T)$. 
\end{remark}

Now, we will prove that the functions $\mathcal{L}_{-}^{\lambda_{il}}\gamma_{il}(x,t)$ and $
\mathcal{L}_{+}^{\beta_{jp}}\theta_{jp}(x,t)$ are well defined. By Lemma \ref{edbf}, it suffices to show that $\gamma_{il},\theta_{jp}\in H_0^{\frac{s+2}{5}}(\R^+)$. Using the expression \eqref{M}, we see that these functions are linear combinations of $$-F_1(0,t)+a_iF_i(0,t), -F_1(0,t)+a_{k+j} G_j(0,t),$$ $$-\partial_x \mathcal{I}_{\frac{1}{5}}F_1(0,t)+b_i\partial_x \mathcal{I}_{\frac{1}{5}}F_i(0,t),$$ $$-\partial_x \mathcal{I}_{\frac{1}{5}}F_1(0,t)+b_{k+j}\partial_x \mathcal{I}_{\frac{1}{5}}G_j(0,t),$$ $$\partial_x^2 \mathcal{I}_{\frac{2}{5}}F_i(0,t)-c_{i1}\partial_x^2 \mathcal{I}_{\frac{2}{5}}G_1(0,t)-\hdots-c_{im}\partial_x^2 \mathcal{I}_{\frac{2}{5}}G_m(0,t),$$ $$-\displaystyle\sum_{i=1}^k b_i^{-1}\partial_x^3 \mathcal{I}_{\frac{3}{5}}F_i(0,t) + \sum_{j=1}^m b_{k+j}^{-1}\partial_x^3 \mathcal{I}_{\frac{3}{5}}G_j(0,t)\quad\text{and}$$ $$-\displaystyle\sum_{i=1}^k a_i^{-1}\partial_x^4 \mathcal{I}_{\frac{4}{5}}F_i(0,t) + \sum_{j=1}^m a_{k+j}^{-1}\partial_x^4 \mathcal{I}_{\frac{4}{5}}G_j(0,t),$$ which we must show are in the appropriate  space. By using Lemmas \ref{grupok}, \ref{edbf}, \ref{duhamelk} and \ref{bilinear1}, we obtain
\begin{equation}
\|F_i(0,t)\|_{H^{\frac{s+2}{5}}(\R^+)}\leq c( \|u_{0i}\|_{H^s(\R^+)}+\|u_{i}\|_{X^{s,b}}^2+\|u_{i}\|_{D^{\alpha}}^2),
\end{equation}
\begin{equation}
\|G_j(0,t)\|_{H^{\frac{s+2}{5}}(\R^+)}\leq c( \|v_{0j}\|_{H^s(\R^+)}+\|v_{j}\|_{X^{s,b}}^2+\|v_{j}\|_{D^{\alpha}}^2).
\end{equation}
If $-\frac{1}{2}<s<\frac{1}{2}$, we have $\frac{3}{10} <\frac{s+2}{5}<\frac12$. Thus Lemma \ref{sobolevh0} implies that $H^{\frac{s+2}{5}}(\R^+)=H_0^{\frac{s+2}{5}}(\R^+)$. It follows that $F_i(0,t),G_j(0,t)\in H_0^{\frac{s+2}{5}}(\R^+)$ for $-\frac{1}{2}<s<\frac{1}{2}.$\\

If $\frac{1}{2}<s<\frac{5}{2}$, then $\frac12<\frac{s+2}{5}<\frac{9}{10}$. Using the compatibility condition \eqref{compatibility1} in Theorem \ref{grandeteorema}, we have that 
\begin{equation*}
F_1(0,0)-a_i F_i(0,0)=u_1(0,0)-a_iu_i(0,0)=u_{0,1}(0)-a_iu_{0,i}(0)=0,
\end{equation*}
and
\begin{equation*}
F_1(0,0)-a_{k+j} G_j(0,0)=u_1(0,0)-a_{k+j}v_j(0,0)=u_{0,1}(0)-a_{k+j}v_{0,j}(0)=0.
\end{equation*}
Then Lemma \ref{alta} implies 
\begin{equation}\label{trace1}
\begin{split}
F_1(0,t)-a_iF_i(0,t)\in H_0^{\frac{s+2}{5}}(\R^+)\quad\text{and}\\
F_1(0,t)-a_{k+j}G_j(0,t)\in H_0^{\frac{s+2}{5}}(\R^+).
\end{split}
\end{equation}

Similarly, using Lemmas \ref{grupok}, \ref{edbf}, \ref{duhamelk} and \ref{bilinear1}, we see that
\begin{equation*}
\|\partial_x F_i(0,t)\|_{H^{\frac{s+1}{5}}(\R^+)}\leq c( \|u_{0,i}\|_{H^s(\R^+)}+\|u_i\|_{X^{s,b}}^2+\|u_i\|_{D^{\alpha}}^2), 
\end{equation*}
\begin{equation*}
\|\partial_xG_j(0,t)\|_{H^{\frac{s+1}{5}}(\R^+)}\leq c( \|v_{0,j}\|_{H^s(\R^+)}+\|v_j\|_{X^{s,b}}^2+\|v_j\|_{D^{\alpha}}^2). 
\end{equation*}

If $-\frac{1}{2}<s<\frac{3}{2}$, we have $\frac{1}{10} <\frac{s+1}{5}<\frac12$. Thus Lemma \ref{sobolevh0} implies that $H^{\frac{s+1}{5}}(\R^+)=H_0^{\frac{s+1}{5}}(\R^+)$. It follows that $\partial_x F_i(0,t),\partial_x G_j(0,t)\in H_0^{\frac{s+1}{5}}(\R^+)$ for $-\frac{1}{2}<s<\frac{3}{2}.$

If $\frac{3}{2}<s<\frac{5}{2}$, then $\frac12<\frac{s+1}{5}<\frac7{10}$. Using the compatibility condition \eqref{compatibility2} in Theorem \ref{grandeteorema}, we have that 
\begin{equation*}
\partial_x F_1(0,0)-b_i \partial_x F_i(0,0)=\partial_x u_1(0,0)-b_i\partial_x u_i(0,0)=\partial_x u_{0,1}(0)-b_i\partial_x u_{0,i}(0)=0,
\end{equation*}
and
\begin{equation*}
\partial_x F_1(0,0)-b_{k+j} \partial_x G_j(0,0)=\partial_x u_1(0,0)-b_{k+j}\partial_x v_j(0,0)=\partial_x u_{0,1}(0)-b_{k+j}\partial_x v_{0,j}(0)=0.
\end{equation*}
Then Lemma \ref{alta} implies $\partial_x F_1(0,t)-b_i\partial_x F_i(0,t),\partial_x F_1(0,t)-b_{k+j}\partial_x G_j(0,t)\in H_0^{\frac{s+1}{5}}(\R^+)$. Therefore, consequence of Lemma \ref{lio-lemaint}, we conclude that
\begin{equation}\label{trace2}
\begin{split}
\partial_x\mathcal{I}_{\frac15} F_1(0,t)-b_i\partial_x\mathcal{I}_{\frac15} F_i(0,t)\in H_0^{\frac{s+2}{5}}(\R^+)\quad\text{and}\\
\partial_x\mathcal{I}_{\frac15} F_1(0,t)-b_{k+j}\partial_x\mathcal{I}_{\frac15} G_j(0,t)\in H_0^{\frac{s+2}{5}}(\R^+).
\end{split}
\end{equation}
Now, again Lemmas  \ref{grupok}, \ref{edbf}, \ref{duhamelk} and \ref{bilinear1} imply
\begin{equation*}
\|\partial^2_x F_i(0,t)\|_{H^{\frac{s}{5}}(\R^+)}\leq c( \|u_{0,i}\|_{H^s(\R^+)}+\|u_i\|_{X^{s,b}}^2+\|u_i\|_{D^{\alpha}}^2), 
\end{equation*}
\begin{equation*}
\|\partial^2_xG_j(0,t)\|_{H^{\frac{s}{5}}(\R^+)}\leq c( \|v_{0,j}\|_{H^s(\R^+)}+\|v_j\|_{X^{s,b}}^2+\|v_j\|_{D^{\alpha}}^2). 
\end{equation*}
Since $-\frac12<s<\frac52$ we have $-\frac{1}{10}<\frac{s}{5}<\frac1{2}$, then Lemma \ref{sobolevh0} gives $\partial^2_xF_i(0,t),\partial^2_xG_j(0,t)\in H_0^{\frac{s}{5}}(\R^+)$. Then using Lemma \ref{lio-lemaint}, we have
\begin{equation}\label{trace3}
    \partial^2_x \mathcal{I}_{\frac25}F_i(0,t),\partial^2_x\mathcal{I}_{\frac25}G_j(0,t)\in H_0^{\frac{s+2}{5}}(\R^+).
\end{equation}

In the same way, we can obtain
\begin{equation}\label{trace4}
    \partial_x^3\mathcal{I}_{\frac{3}{5}}F_i(0,t), \partial_x^3\mathcal{I}_{\frac{3}{5}}G_j(0,t), \partial_x^4\mathcal{I}_{\frac{4}{5}}F_i(0,t), \partial_x^4\mathcal{I}_{\frac{4}{5}}G_j(0,t) \in H_0^{\frac{s+2}{5}}(\R^+).
\end{equation}
Thus, \eqref{trace1}, \eqref{trace2}, \eqref{trace3} and \eqref{trace4} imply that the functions ${\mathcal{V}_{-}^{\lambda_{il}}\gamma_{il}(x,t)}$ {and} $
{\mathcal{V}_{+}^{\beta_{jp}}\theta_{jp}(x,t)}$ are well defined.

\section{End of the proof: Obtaining a fixed point of $\Lambda$}
Using Lemmas  \ref{lio}, \ref{grupok}, \ref{edbf}, \ref{duhamelk} and \ref{bilinear1}, we obtain
\begin{equation}
\begin{split}
\|\Lambda (\tilde{u}_{1},\hdots,\tilde{u}_{k},\tilde{v}_{1},\hdots,\tilde{v}_{m})-\Lambda (\tilde{\tilde{u}}_{1},\hdots,\tilde{\tilde{u}}_{k},\tilde{\tilde{v}}_{1},\hdots,\tilde{\tilde{v}}_{m})\|_{Z_{(s,b,\alpha)}}\\
\leq c(\|\Lambda (\tilde{u}_{1},\hdots,\tilde{u}_{k},\tilde{v}_{1},\hdots,\tilde{v}_{m})\|_{Z_{(s,b,\alpha)}}+\|\Lambda (\tilde{\tilde{u}}_{1},\hdots,\tilde{\tilde{u}}_{k},\tilde{\tilde{v}}_{1},\hdots,\tilde{\tilde{v}}_{m})\|_{Z_{(s,b,\alpha)}})\\ \cdot\|\Lambda (\tilde{u}_{1},\hdots,\tilde{u}_{k},\tilde{v}_{1},\hdots,\tilde{v}_{m})-\Lambda (\tilde{\tilde{u}}_{1},\hdots,\tilde{\tilde{u}}_{k},\tilde{\tilde{v}}_{1},\hdots,\tilde{\tilde{v}}_{m})\|_{Z_{(s,b,\alpha)}}
\end{split}
\end{equation}
and 
\begin{equation}
\begin{split}
\|\Lambda ({u}_{1},\hdots,{u}_{k},{v}_{1},\hdots,{v}_{m})\|_{Z}&\leq c \left( \sum_{i=1}^k\|u_{0,i}\|_{H^s(\R^+)}+\sum_{j=1}^m\| v_{0,j}\|_{H^s(\R^+)}\right.\\
&\left.+\sum_{i=1}^k(\|u_i\|_{X^{s,b}}^2+\|u_i\|_{D^{\alpha}}^2)+\sum_{j=1}^m(\|v_j\|_{X^{s,b}}^2+\|v_j\|_{D^{\alpha}}^2)\right) .
\end{split}
\end{equation}

By taking $\sum_{i=1}^{k}\|u_{0,i}\|_{H^s(\R^+)}+\sum_{j=1}^{m}\|v_{0,j}\|_{H^s(\R^+)}<\delta$ for $\delta>0$ suitable small, we obtain a fixed point $\Lambda (\overline{u}_{1},\hdots,\overline{u}_{k},\overline{v}_{1},\hdots,\overline{v}_{m})=(\overline{u}_{1},\hdots,\overline{u}_{k},\overline{v}_{1},\hdots,\overline{v}_{m})$ in a ball $$B=\{({u}_{1},\hdots,{u}_{k},{v}_{1},\hdots,{v}_{m})\in Z; \|({u}_{1},\hdots,{u}_{k},{v}_{1},\hdots,{v}_{m})\|_Z\leq 2c\delta\}.$$ It follows that the restriction
\begin{equation}
({u}_{1},\hdots,{u}_{k},{v}_{1},\hdots,{v}_{m})=(\overline{u_1}\big|_{\R^-\times (0,1)},\hdots,\overline{u_k}\big|_{\R^-\times (0,1)},\overline{v_1}\big|_{\R^+\times (0,1)},\hdots, \overline{v_m}\big|_{\R^+\times (0,1)})
\end{equation}
solves the Cauchy problem \eqref{Kawahara} with boundary conditions \eqref{condfront}, on the space $\mathcal Z_{s,b,\alpha}$, in the sense of distributions.

Finally, the existence of solutions for any data in $H^s(\mathcal{G})$ follows by the standard scaling argument, the reader can consult \cite{MC}. 

\medskip

\section{Applications of Theorem \ref{grandeteorema}}
In this section, we give applications of the general criteria described in Theorem \ref{grandeteorema} in some situations, including balanced and unbalanced star graphs.

\begin{example}[$\mathcal{Y}$-junction] Consider the $\mathcal{Y}$-junction, i.e, the case $k=1$ and $m=2$, with boundary conditions are given by $a_1=a_{2}=a_3=1$, $b_1=b_{2}=b_3=1$ and $C=\left[\frac{\sqrt{2}}{2}\quad \frac{\sqrt{2}}{2}\right]$. For each $s\in (-0.40,2.40)\setminus\{0.5,1.5\}$, there exists a time $T>0$ and a solution of the Cauchy problem \eqref{Kawahara}-\eqref{condfront}  satisfying the compatibility conditions \eqref{compatibility1} when $s>0.5$ and \eqref{compatibility2} $s>1.5$.
\end{example}
\begin{figure}[H]
    \centering
    \begin{tikzpicture}[line cap=round]
        \draw[dashed, black, very thick] (0,-3) -- (0,-2);
        \draw[thick, black, very thick] (0,-2) -- (0,0);
        \draw[black, ->, thick] (0,-2.3) -- (0,-2.2);
        \draw[black, ->, thick] (0,-1.5) -- (0,-1.2);
        \draw[thick, red, very thick] (0,0) -- (135:3);
        \draw[dashed, red, very thick] (135:3) -- (135:4);
        \draw[red, ->, thick] (135:1) -- (135:1.5);
        \draw[red, ->, thick] (135:3.5) -- (135:3.75);
        \draw[thick, red, very thick] (0,0) -- (45:3);
        \draw[dashed, red, very thick] (45:3) -- (45:4);
        \draw[red, ->, thick] (45:1) -- (45:1.5);
        \draw[red, ->, thick] (45:3.5) -- (45:3.75);
        \node at (-3.5,1.5) {$( 0, +\infty)$};
        \node at (3.5,1.5) {$(0, +\infty)$};
        \node at (1.5,-1.7) {$(-\infty,0)$};
    \end{tikzpicture}
    \caption{$\mathcal{Y}$-junction}
    \label{fig-y-junction}
\end{figure}
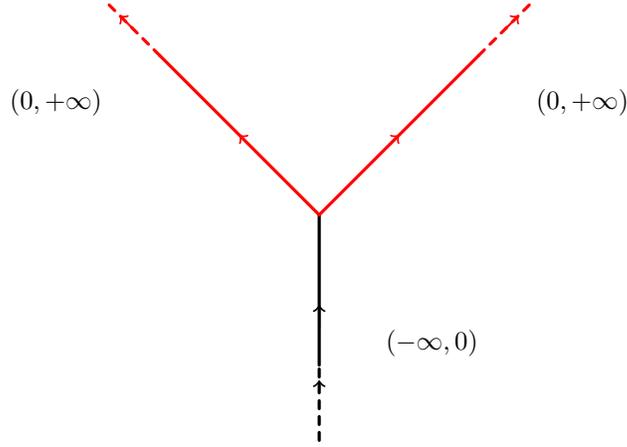
To verify this, we calculate the determinant of
$\mathbf{M}(\boldsymbol\lambda, \boldsymbol\beta)$ using the code in the appendix. Those calculations were made for the following:
\begin{itemize}
    \item $\lambda_{1,1}=0.48$, $\lambda_{1,2}=0.32$, $\lambda_{1,3}=0.16$, $\beta_{1,1}=\beta_{1,2}=0.48 $ and $\beta_{2,1}=\beta_{2,2}=0.23 $, resulting in $\det(\mathbf{M}(\boldsymbol\lambda, \boldsymbol\beta))=-9.2205 \cdot\left[\sin{\left(\frac{\pi}{10}\right)}\right]^{7}$;
    \item $\lambda_{1,1}=0.49$, $\lambda_{1,2}=0.45$, $\lambda_{1,3}=0.40$, $\beta_{1,1}=\beta_{1,2}=0.49 $ and $\beta_{2,1}=\beta_{2,2}=0.40$, resulting in $\det(\mathbf{M}(\boldsymbol\lambda, \boldsymbol\beta))=-0.0740 \cdot\left[\sin{\left(\frac{\pi}{10}\right)}\right]^{7}$;
    \item $\lambda_{1,1}=0.10$, $\lambda_{1,2}=0.05$, $\lambda_{1,3}=0.01$, $\beta_{1,1}=\beta_{1,2}=0.10 $ and $\beta_{2,1}=\beta_{2,2}=0.01$, resulting in $\det(\mathbf{M}(\boldsymbol\lambda, \boldsymbol\beta))=-0.0836 \cdot\left[\sin{\left(\frac{\pi}{10}\right)}\right]^{7}$.
\end{itemize}
The first calculation validates the result for $s\in(-0.02,2.16)\setminus\{0.5,1.5\}$; the second one, for $s\in(-0.01,2.40)\setminus\{0.5,1.5\}$; and the third, for $s\in(-0.40,2.01)\setminus\{0.5,1.5\}$.

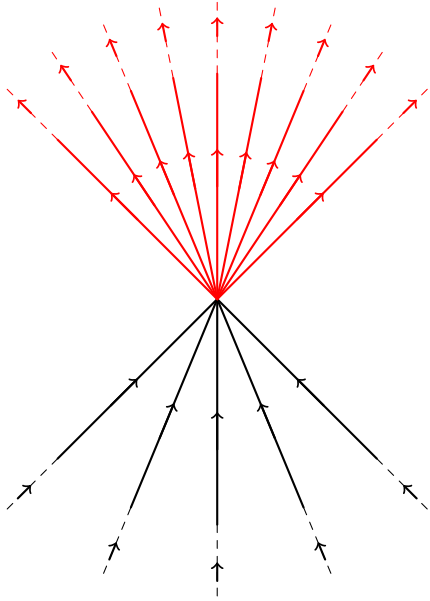
\begin{figure}[H]
    \centering
\begin{tikzpicture}[line cap=round]
\def\numLeftEdges{5}
\def\numRightEdges{9}
\def\angleSpread{90}
\foreach \i in {1,...,\numLeftEdges} {
  \draw[thick, black] (0,0) -- ({270 - \angleSpread/2 + (\i-1)*\angleSpread/(\numLeftEdges-1)}:3); 
  \draw[dashed, black] ({270 - \angleSpread/2 + (\i-1)*\angleSpread/(\numLeftEdges-1)}:3) -- ({270 - \angleSpread/2 + (\i-1)*\angleSpread/(\numLeftEdges-1)}:4) node[left, black] {};
  \draw[black, <-, thick] ({270 - \angleSpread/2 + (\i-1)*\angleSpread/(\numLeftEdges-1)}:1.5) -- ({270 - \angleSpread/2 + (\i-1)*\angleSpread/(\numLeftEdges-1)}:2); 
  \draw[black, <-, thick] ({270 - \angleSpread/2 + (\i-1)*\angleSpread/(\numLeftEdges-1)}:3.5) -- ({270 - \angleSpread/2 + (\i-1)*\angleSpread/(\numLeftEdges-1)}:3.75);
}
\foreach \i in {1,...,\numRightEdges} {
  \draw[thick, red] (0,0) -- ({90 + \angleSpread/2 - (\i-1)*\angleSpread/(\numRightEdges-1)}:3); 
  \draw[dashed, red] ({90 + \angleSpread/2 - (\i-1)*\angleSpread/(\numRightEdges-1)}:3) -- ({90 + \angleSpread/2 - (\i-1)*\angleSpread/(\numRightEdges-1)}:4) node[right, red] {};
  \draw[red, ->, thick] ({90 + \angleSpread/2 - (\i-1)*\angleSpread/(\numRightEdges-1)}:1.5) -- ({90 + \angleSpread/2 - (\i-1)*\angleSpread/(\numRightEdges-1)}:2);
  \draw[red, ->, thick] ({90 + \angleSpread/2 - (\i-1)*\angleSpread/(\numRightEdges-1)}:3.5) -- ({90 + \angleSpread/2 - (\i-1)*\angleSpread/(\numRightEdges-1)}:3.75);
}
\end{tikzpicture}
\caption{A star graph with five $(-\infty,0)$ edges (the black ones) and eight $(0,\infty)$ edges (the red ones)}
\end{figure}

\begin{example}[Others star graphs]\label{others}
We can also take many other examples in which the Cauchy problem \eqref{Kawahara}-\eqref{condfront} has a solution. For example, for $k,m\leq10$, $a_1=a_2=\hdots=a_{k+m}=1$, $b_1=b_2=\hdots=b_{k+m}=1$, $C= \left[ (km)^{-1/2}\right]_{k\times m}$ and the same hypothesis of Theorem \ref{grandeteorema}, we get solutions as in that theorem for $s\in (-0.02,2.16)\setminus\{0.5,1.5\}$. Table \ref{tabeladassimulacoes} shows by computer simulations that the matrix $\mathbf{M}(\boldsymbol\lambda, \boldsymbol\beta)$ is invertible.
\end{example}

\begin{table}[H]
\begin{tabular}{lrrrrrrrrrr}
 & 1 & 2 & 3 & 4 & 5 & 6 & 7 & 8 & 9 & 10 \cr
   1 & -3.627 & -9.221 & 23.147 & 55.247 & -123.968 & -263.127 & 533.350 & 1041.356 & -1972.164 & -3642.499 \cr
    2 & 2.212 & -4.309 & -8.927 & 18.925 & 39.677 & -81.049 & -160.772 & 310.135 & 583.530 & -1074.316 \cr
    3 & 1.286 &  2.099 & -3.739 & -7.074 & 13.742 & 26.746 & -51.488 & -97.501 & 181.379 & 331.641 \cr
   4 & -0.704 &  1.032 &  1.650 & -2.828 & -5.079 &  9.343 & 17.304 & -31.935 & -58.404 & 105.596 \cr
   5 & -0.361 & -0.498 &  0.742 &  1.179 & -1.975 & -3.429 &  6.077 & 10.863 & -19.421 & -34.562 \cr
    6 & 0.175 & -0.234 & -0.333 &  0.502 &  0.794 & -1.308 & -2.218 &  3.827 &  6.660 & -11.617 \cr
    7 & 0.081 &  0.107 & -0.148 & -0.214 &  0.325 &  0.513 & -0.835 & -1.390 &  2.350 &  4.005 \cr
   8 & -0.036 &  0.047 &  0.064 & -0.091 & -0.134 &  0.204 &  0.321 & -0.517 & -0.850 &  1.414 \cr
   9 & -0.016 & -0.020 &  0.027 &  0.038 & -0.055 & -0.082 &  0.125 &  0.196 & -0.314 & -0.509 \cr
    10 & 0.007 & -0.009 & -0.012 &  0.016 &  0.023 & -0.033 & -0.049 &  0.075 &  0.118 & -0.187
\end{tabular}
\caption{Values of $det(\mathbf{M}(\boldsymbol\lambda, \boldsymbol\beta))\cdot\left[\sin{\left(\frac{\pi}{10}\right)}\right]^{-(3k+2m)}$ for $\lambda_{i1}=0.48,\lambda_{i2}=0.32,\lambda_{i3}=0.16,\beta_{j1}=0.48,\beta_{j2}=0.23,$ and $a_i,a_{k+j},b_i,b_{k+j}$ and $C$ as in Example \ref{others}. Note that all of them are non-null. The row represents the value of $k$ and the column the value of $m$.}
\label{tabeladassimulacoes}
\end{table}

\begin{remark}
In the above examples, one can expand the range in which $s$ can be taken. To do this, it suffices to look for different parameters $\lambda_{il}$ and $\beta_{jp}$ for which the matrix $\mathbf{M}(\boldsymbol\lambda, \boldsymbol\beta)$ is invertible.
\end{remark}

\section*{Acknowledgments}
 M. Cavalcante wishes to thank the support of CNPq, 
grant \# 310271/
2021-5, and the support of the Funda\c c\~ao de Amparo \`a Pesquisa do Estado de Alagoas
- FAPEAL, Brazil, grant \# E:60030.0000000161/2022. The authors thank the CAPES/Cofecub, grant \# 88887.879175/2023-00.
\section*{Statements and Declaration}

\textbf{Conflicts of interest} All authors declare that they have no conflicts of interest.
\medskip
\appendix
\section{Code for computation of the determinant}
A MATLAB/Octave function to calculate the determinant of the matrix $M(\lambda, \beta)$ follows. In the input, $X$ means the $k\times 3$-dimensional matrix with general term $x_{il}=\lambda_{il}$ and $Y$ means the $m\times 2$-dimensional matrix, the general term being $y_{jp}=\beta_{jp}$.\\

Furthermore, the meaning of $a(u,v),b(u,v),\hdots,e(u,v)$ is explained in the following table:

\begin{table}[H]
\begin{tabular}{ll}
 variable on code & variable on theorem\\
 $a(u,v)$ & $\rho^{(u-1)}_{v1}$\\[4pt]
 $b(u,v)$ & $\rho^{(u-1)}_{v2}$\\[4pt]
 $c(u,v)$ & $\rho^{(u-1)}_{v3}$\\[4pt]
 $d(u,v)$ & $\iota^{(u-1)}_{v1}$\\[4pt]
 $e(u,v)$ & $\iota^{(u-1)}_{v2}$\\[4pt]
\end{tabular}
\end{table}
\medskip

\begin{lstlisting}[language=MATLAB] 
% code for calculating the determinant
function retval = detergeral(X,Y,A,B,C)
k = size(X,1);
m = size(Y,1);


%building the matrix
i=0.0000 + 1.0000i;
for u=1:5
for v = 1:m
	d(u,v)=4*sin(pi/10)*cos(-Y(v,1)*pi/5+u*pi/5);
	e(u,v)=4*sin(pi/10)*cos(-Y(v,2)*pi/5+u*pi/5);
    end
  for v = 1:k
	a(u,v)=2*sin(pi/10)*[2*cos(-2*X(v,1)*pi/5+2*u*pi/5)+1];
	b(u,v)=2*sin(pi/10)*[2*cos(-2*X(v,2)*pi/5+2*u*pi/5)+1];
        c(u,v)=2*sin(pi/10)*[2*cos(-2*X(v,3)*pi/5+2*u*pi/5)+1];
  end
  end


M=zeros(3*k+2*m);

% matrices M11 and M21 
M(1:k+m-1,1)=a(1,1)*ones(k+m-1,1);
M(1:k+m-1,2)=b(1,1)*ones(k+m-1,1);
M(1:k+m-1,3)=b(1,1)*ones(k+m-1,1);
for v = 2:k
	M(v-1,3*v-2) = -A(v)*a(1,v);
	end
for v = 2:k
	M(v-1,3*v-1) = -A(v)*b(1,v);
        end
for v = 2:k
	M(v-1,3*v) = -A(v)*c(1,v);
        end
% matrix M22
for v = 1:m
	M(k-1+v,3*k+2*v-1) = -A(k+v)*d(1,v);
	end
for v = 1:m
	M(k-1+v,3*k+2*v) = -A(k+v)*e(1,v);
        end
% matrices M31 and M41 
M(k+m:2*k+m-1,1)=a(2,1)*ones(k+m-1,1);
M(k+m:2*k+m-1,2)=b(2,1)*ones(k+m-1,1);
M(k+m:2*k+m-1,3)=c(2,1)*ones(k+m-1,1);
for v = 2:k
	M(k+m+v-2,3*v-2) = -B(v)*a(2,v);
	end
for v = 2:k
	M(k+m+v-2,3*v-1) = -B(v)*b(2,v);
        end
for v = 2:k
	M(k+m+v-2,3*v) = -B(v)*c(2,v);
	end
% matrix M42
for v = 1:m
	M(2*k+m-2+v,3*k+2*v-1) = -B(k+v)*d(2,v);
	end
for v = 1:m
	M(2*k+m-2+v,3*k+2*v) = -B(k+v)*e(2,v);
        end
% matrix M51
for v = 1:k
	M(2*k+2*m-2+v,3*v-2) = -a(3,v);
	end
for v = 1:k
	M(2*k+2*m-2+v,3*v-1) = -b(3,v);
	end
for v = 1:k
	M(2*k+2*m-2+v,3*v) = -c(3,v);
	end
% matrix M52
for v = 1:m
	M(2*k+2*m-1:3*k+2*m-2,3*k+2*v-1) = d(3,v)*C(1:k,v);
        end
for v = 1:m
	M(2*k+2*m-1:3*k+2*m-2,3*k+2*v) = e(3,v)*C(1:k,v);
        end
% matrix M61
for v = 1:k
	M(3*k+2*m-1,3*v-2) = ((B(v))^(-1))*a(4,v);
	end
for v = 1:k
	M(3*k+2*m-1,3*v-1) = ((B(v))^(-1))*b(4,v);
	end
for v = 1:k
	M(3*k+2*m-1,3*v) = ((B(v))^(-1))*c(4,v);
        end
% matrix M62
for v = 1:m
	M(3*k+2*m-1,3*k+2*v-1) = -((B(k+v))^(-1))*d(4,v);
	end
for v = 1:m
	M(3*k+2*m-1,3*k+2*v) = -((B(k+v))^(-1))*e(4,v);
	end
% matrix M71
for v = 1:k
	M(3*k+2*m,3*v-2) = ((A(v))^(-1))*a(5,v);
	end
for v = 1:k
	M(3*k+2*m,3*v-1) = ((A(v))^(-1))*b(5,v);
	end
for v = 1:k
	M(3*k+2*m,3*v) = ((A(v))^(-1))*c(5,v);
        end
% matrix M72
for v = 1:m
	M(3*k+2*m,3*k+2*v-1) = -((A(k+v))^(-1))*d(5,v);
	end
for v = 1:m
	M(3*k+2*m,3*k+2*v) = -((A(k+v))^(-1))*e(5,v);
	end
size(M);
retval=det(M);

end
\end{lstlisting}

\end{document}